
\magnification=1200
\hfuzz=3pt
\overfullrule=0mm

\hsize=125mm
\hoffset=4mm


\font\tensymb=msam9
\font\fivesymb=msam5 at 5pt
\font\sevensymb=msam7  at 7pt
\newfam\symbfam
\scriptscriptfont\symbfam=\fivesymb
\textfont\symbfam=\tensymb
\scriptfont\symbfam=\sevensymb

\font\titlefont=cmbx10 at 15pt


\font\refttfont=cmtt10 at 9pt

\font\sc=cmcsc10 \rm


\def\Hom{{\rm Hom}}
\def\Inf{{\rm Inf}}
\def\Ker{{\rm Ker}}
\def\Im{{\rm Im}}
\def\id{{\rm id}}
\def\free{{\rm free}}
\def\univ{{\rm univ}}
\def\ot{\otimes}

\def\ZZ{{\bf Z}}

\def\cE{{\cal E}}
\def\cF{{\cal F}}

\def\and{\quad\hbox{and}\quad}

\def\Pr{\noindent {\sc Proof.--- }}

\def\cqfd{ {\sevensymb {\char 3}}}

\def\hfl#1#2{\smash{\mathop{\hbox to 6mm{\rightarrowfill}}
\limits^{\scriptstyle#1}_{\scriptstyle#2}}}

\def\vfl#1#2{\llap{$\scriptstyle #1$}\left\downarrow
\vbox to 3mm{}\right.\rlap{$\scriptstyle #2$}}


\null

\vskip 15pt

\centerline{\titlefont Norm formulas for finite groups and induction}
\medskip
\centerline{\titlefont from elementary abelian subgroups}

\vskip 30pt

\centerline{\sc Eli Aljadeff and Christian Kassel}

\bigskip\bigskip
\noindent
{\sc Abstract.}
{\it It is known that the norm map~$N_G$ for a finite
group~$G$ acting on a ring~$R$ is surjective if and only if
for every elementary abelian subgroup $E$ of~$G$
the norm map~$N_E$ for~$E$ is surjective. Equivalently, there exists
an element $x_G\in R$ with $N_G(x_G) = 1$ if and only
for every elementary abelian subgroup~$E$ there exists an
element~$x_E\in R$ such that $N_E(x_E) = 1$.
When the ring $R$ is noncommutative, it is an open problem to find an
explicit formula
for $x_G$ in terms of the elements~$x_E$.
In this paper we present a method to solve this problem for an arbitrary
group~$G$ and an arbitrary group action on a ring.
Using this method, we obtain a complete
solution of the problem for the quaternion and the dihedral $2$-groups,
and for a group of order~27.
We also show how to reduce the problem to the class
of almost extraspecial $p$-groups.
}



\vskip 15pt
\bigskip

\noindent

Let $G$ be a finite group acting by ring automorphisms 
on an arbitrary (not necessarily commutative) ring~$R$ with unit.
For any subgroup $U$ of~$G$ the {\it norm map}
$N_U : R \to R^U$ is defined for all $x\in R$ by
$$N_U(x) = \sum_{g\in U}\, g(x).$$
Here $g(x)$ denotes the value in~$R$ of the action of $g$ on~$x$
and $R^U$ the subring of $U$-invariant elements in~$R$.

The question of the surjectivity of the map $N_G$ onto~$R^G$ has
well-known interpretations in topics such as Galois theory,
algebraic number fields and, most importantly for this paper,
integral group representations.

In [2, Theorem~1] the first-named author and Ginosar proved that 
$N_G$ is surjective onto~$R^G$ if and only if $N_E$ is surjective onto~$R^E$ 
for every elementary abelian subgroup $E$ of~$G$.
This generalizes Chouinard's theorem [7] that asserts that
a $\ZZ[G]$-module is projective if and only if 
for every elementary abelian subgroup $E$ of~$G$
it is projective as a $\ZZ[E]$-module.

The map $N_U$ being $R^U$-linear, it is surjective onto~$R^U$
if and only if the unit $1$ of~$R$ belongs to the image of~$N_U$.
We can therefore rephrase Aljadeff and Ginosar's result as follows: 
there is an element~$x_G \in R$ such that $N_G(x_G) = 1$ if and only if
there is an element~$x_E\in R$ such that $N_E(x_E) = 1$
for every elementary abelian subgroup~$E$ of~$G$.
Using this statement, Shelah observed (see [2,~Proposition~6])
that there exist formulas expressing~$x_G$
polynomially in terms of the elements~$x_E$ and the elements of~$G$.

The aim of this paper is to find such a formula for any given finite group~$G$.
We would like to point out that such a formula should be defined over $\ZZ$
and be independent of the ring on which the group acts. 
In this way it has a universal character though it may not be unique.
With such a formula it is possible to construct a projective $G$-basis
for any finitely generated projective $\ZZ[G]$-module~$M$
out of given projective $E$-bases of~$M$, one for each
elementary abelian subgroup $E$ of~$G$.

We assume that $G$ is not an elementary abelian group
(otherwise the problem is trivial).
The smallest group that is not elementary abelian
is the cyclic group $C_4$ of order~$4$.
This case was solved by P\'eter P.~P\'alfy who found the formula
$$x_G = x_E\sigma(x_E) + x_E\sigma(x_E)x_E  - x_E^2\sigma(x_E), \eqno (0.1)$$
where $\sigma$ denotes a generator of~$C_4$.
The meaning of this formula is the following: if $C_4$ acts
on some ring~$R$ and there is an element $x_E \in R$
such that $x_E + \sigma^2(x_E) = 1$ (i.e., $N_E(x_E) = 1$
for the (unique) elementary abelian subgroup~$E$ of~$G$),
then for $x_G\in R$ given by~(0.1) we have
$$N_G(x_G) = x_G + \sigma(x_G) + \sigma^2(x_G) + \sigma^3(x_G) = 1.$$
This was the first formula of this kind to be found.

The next step is due to the authors:  they showed in~[3]
how to obtain formulas for {\it all abelian groups}.
Such formulas can be complicated even in simple cases.
For instance, if $G = C_9$ is a cyclic group of order~$9$
(with generator~$\sigma$)
and $x_E\in R$ is an element of norm one for the
subgroup of order~$3$, then
$$\eqalign{
x_G = &
- x_E^2 + 2 \sigma(x_E)x_E - \sigma^3(x_E)x_E + \sigma^4(x_E)x_E \cr
&{} + x_E\sigma^3(x_E)x_E + x_E\sigma^4(x_E)x_E + x_E\sigma^5(x_E)x_E \cr
&{} + x_E\sigma^6(x_E)x_E + x_E\sigma^7(x_E)x_E + x_E\sigma^8(x_E)x_E\cr
&{} - \sigma(x_E)\sigma^4(x_E)x_E - \sigma(x_E)\sigma^5(x_E)x_E
- \sigma(x_E)\sigma^6(x_E)x_E \cr
&{} - \sigma(x_E)\sigma^7(x_E)x_E -
\sigma(x_E)\sigma^8(x_E)x_E - \sigma(x_E) x_E^2 \cr
& + \sigma^3(x_E)\sigma^6(x_E)x_E + \sigma^3(x_E)\sigma^7(x_E)x_E +
\sigma^3(x_E)\sigma^8(x_E)x_E \cr
& - \sigma^4(x_E)\sigma^7(x_E)x_E -  \sigma^4(x_E)\sigma^8(x_E)x_E
- \sigma^4(x_E)x_E^2\cr
} \eqno (0.2)$$
is a formula (with 22 monomials) for an element of norm one for~$G$.
It should also be noted that the first-named author obtained formulas
for arbitrary groups acting on {\it commutative} rings (see~[1]).
\goodbreak

We are thus left with the case of {\it nonabelian groups acting on noncommutative rings}
(noncommutative rings are important for us because we want to be able to apply the formulas
to rings of endomorphisms).
In this paper we present a general method that allows to solve the problem
for arbitrary groups and arbitrary ring actions. The idea is to
translate the task of finding formulas for norm one elements $x_G$
for a group~$G$ acting on a ring~$R$
into a more accessible system of equations in~$R$, where the variables
are symbols~$b(\sigma)$, one for each generator $\sigma$ in a presentation
of~$G$. The equations are obtained from the relations in the presentation.
Such a system has solutions, and we will show how each solution yields
a formula for~$x_G$.
We illustrate this method by solving this system of equations
for the quaternion and dihedral families of $2$-groups, and a group of
order~27.
This provides the first examples of norm one formulas for nonabelian groups
acting on noncommutative rings.

We also show that the problem for a general group $G$ can be reduced
to a smaller class of groups, namely the class of
extraspecial and almost extraspecial $p$-groups that are subquotients of~$G$.
For instance, in order to solve the problem
for the quaternion and dihedral groups mentioned above, it is sufficient
to solve it for the quaternion group $Q_8$ and the dihedral group~$D_8$
of order~$8$.

The paper is organized as follows.
In Section~1 we explain precisely what we mean by a formula for a group~$G$
and we introduce a ring that is universal for the situation under consideration.

In Section~2 we present three reductions,
first a straightforward one to $p$-groups, then
a reduction to extraspecial and almost extraspecial $p$-groups.
Finally we show how to solve the problem for the product of two $p$-groups
once we have solutions for each of them.

After some cohomological preliminaries in Section~3
we explain our method to solve the problem for an arbitrary $p$-group
in Section~4.
This involves solving the above-mentioned system of equations and a further
problem that we solve in Section~5.
More precisely, assuming we are given an element $x\in R$ such that $N_U(x) = 1$,
we will show in Section~5 how to express any $1$-cocycle $\beta : U \to R$
explicitly as a $1$-coboundary, i.e., how to find an explicit formula
for an element $w\in R$ such that $\beta(g)  = g(w) - w$ for all~$g\in U$.
Once we have solved the above-mentioned system of equations and we have
a formula for~$w$, we obtain a complete explicit solution to the problem of
finding a formula for a norm one element for~$G$.

In Sections~6--8 we apply our method to two
important families of nonabelian groups, namely
the quaternion groups $Q_{2^n}$ and the dihedral groups~$D_{2^n}$,
and to a group of order~$27$.

All groups considered in this paper are finite, and all rings have units.
We denote a cyclic group of order $n$ by~$C_n$.

\bigskip\goodbreak
\noindent
{\bf 1.~Formulas for a group}
\medskip

For any group~$G$
we denote $\cE_G$ the set of elementary abelian subgroups of~$G$.
Recall that a group $E$ is elementary abelian if it is isomorphic
to~$C_p^r$ for some prime number $p$ and some integer $r\geq 1$.
Clearly, $\cE_H \subset \cE_G$ if $H$ is a subgroup of~$G$.

\medskip
\noindent
{\sc Definition~1.1}.---
{\it
A formula for a finite group~$G$ is a
polynomial $\Phi_G$ in noncommuting variables $g(x_E)$, 
where $g\in G$ and $E \in \cE_G$, and with coefficients in~$\ZZ$, 
satisfying the following condition:
whenever $G$ acts by automorphisms on a ring~$R$ and $(x_E)_{E\in \cE_G}$
is a family of elements of~$R$ such that $N_E(x_E) = 1$ for all $E \in \cE_G$,
then the element $x_G \in R$ obtained by replacing in~$\Phi_G$ 
each variable $g(x_E)$ by the value of the action
of~$g$ on the element $x_E \in R$ satisfies $N_G(x_G) = 1$.
}
\goodbreak\medskip

In order to clarify Definition~1.1, 
we consider the free noncommutative ring 
$$R_{\free}(G) = \ZZ\, \langle\, g(X_E) \, |\; g\in G,\, E \in \cE_G \,\rangle$$
generated by symbols of the form
$g(X_E)$, where $g$ runs over all elements of~$G$ and 
$E$ runs over all elements of~$\cE_G$.
To simplify notation, we set $e(X_E) = X_E$ when $e$ is the neutral element of~$G$.
The group $G$ acts by ring automorphisms on~$R_{\free}(G)$ as follows: 
if $h\in G$ and $g(X_E)$ is a generator of $R_{\free}(G)$, then
$$h(g(X_E)) = (hg)(X_E)$$
for $g$, $h\in G$, and $E\in \cE_G$.
Let $I$ be the two-sided ideal of $R_{\free}(G)$ generated by all elements of the form
$$\sum_{h\in E}\, (gh)(X_E) - 1$$
for any $g\in G$ and any $E\in \cE_G$.
The ideal is preserved by the $G$-action on~$R_{\free}(G)$.

Let  $R_{\univ}(G)$ be the quotient ring
$$R_{\univ}(G) = R_{\free}(G)/I$$
with the induced $G$-action.
By definition of $R_{\univ}(G)$, for any $E\in \cE_G$ we have
$$N_E(X_E) = \sum_{h\in E}\, h(X_E) = 1.\eqno (1.1)$$

\medskip\goodbreak
\noindent
{\sc Proposition~1.2}.---
{\it Any element $\Phi_G \in R_{\univ}(G)$ such that $N_G(\Phi_G) = 1$
is a formula for the group $G$.
}

\medskip
\noindent
\Pr
First observe that $\Phi_G$ is a polynomial with integer coefficients in noncommutative
variables~$g(X_E)$ indexed by $G\times \cE_G$.
Suppose that $G$ acts on a ring $R$ and that $(x_E)_{E\in \cE_G}$ is a family of elements
of~$R$ such that $N_E(x_E) = 1$ for all $E\in \cE_G$.
Set $f(g(X_E)) = g(x_E)$ for all $g\in G$ and $E\in \cE_G$. Since
for all $g\in G$ and $E\in \cE_G$,
$$\sum_{h\in E}\, (gh)(x_E) = gN_E(x_E) = g(1) = 1$$
in~$R$, 
there is a unique homomorphism of $G$-rings $f: R_{\univ}(G) \to R$
such that $f(g(X_E)) = g(x_E)$ for all $g\in G$ and~$E\in \cE_G$.
Set $x_G = f(\Phi_G)$: this is the element of~$R$ obtained by replacing 
each variable $g(X_E)$ in~$\Phi_G$ by the value of the action
of~$g$ on~$x_E \in R$. We have
$$N_G(x_G) = N_G(f(\Phi_G)) = f(N_G(\Phi_G)) = f(1) = 1. \eqno \hbox{\cqfd}$$
\goodbreak

The proof above also shows that $R_{\univ}(G)$ is the universal $G$-ring 
with a family of elements $(X_E)_{E\in \cE_G}$ such that $N_E(X_E) = 1$.

As we have already pointed out in the introduction, there is a formula for every
finite group~$G$. Let us give a quick proof of this fact using the ring~$R_{\univ}(G)$:
indeed by~(1.1), we have $N_E(X_E) = 1$ 
for every elementary abelian subgroup $E$ of~$G$. 
Therefore by [2, Theorem~1] there exists $\Phi_G \in R_{\univ}(G)$ such that $N_G(\Phi_G) = 1$.
By Proposition~1.2 this is a formula for~$G$.
Since finding such a formula for $G$ amounts to constructing a norm-one element
in~$R_{\univ}(G)$, we see that the problem is a pure group-theoretical question.

The right-hand sides of (0.1) and (0.2) provide formulas for the cyclic
groups $C_4$ and $C_9$, respectively.
If $E$ is an elementary abelian group, then $\Phi_E = x_E$ is clearly a formula for~$E$. 

In our search for formulas for a group~$G$,
the following rephrasing of~[2, Theorem~1] will be useful.
To state it, we need the following notation:
$\cE_G^{\rm max}$ denotes the set of maximal elements of~$\cE_G$ with
respect to
inclusion, and $\cE_G^0$ a subset of $\cE_G^{\rm max}$ such that any
element of $\cE_G^{\rm max}$ is conjugated in~$G$ to exactly one
element of~$\cE_G^0$.

\medskip\goodbreak
\noindent
{\sc Proposition~1.3}.---
{\it Let $G$ be a finite group acting on a ring~$R$ by ring automorphisms.
We assume that $G$ is not elementary abelian.
Then the following statements are equivalent:

1) There exists $x_G\in R$ such that $N_G(x_G) = 1$.

2) For each proper subgroup $U$ of~$G$ there exists $x_U\in R$ with
$N_U(x_U) = 1$.

3) For each $E\in \cE_G$ there exists $x_E\in R$ such that $N_E(x_E) = 1$.

4) For each $E\in \cE_G^{\rm max}$ there exists $x_E\in R$ such that
$N_E(x_E) = 1$.

5) For each $E\in \cE_G^0$ there exists $x_E\in R$ such that $N_E(x_E) = 1$.
}
\medskip

\Pr
1) $\Rightarrow$ 2): Let $\{g_1, \ldots, g_r\}$ be a set of representatives for the
right cosets of $U$ in~$G$. Define
$$x_U = g_1(x_G) + \cdots + g_r(x_G) \in R. \eqno (1.2)$$
Then
$$N_U(x_U) = \sum_{u\in U}\sum_{i=1}^r (ug_i) (x_G) = N_G(x_G) = 1.$$

2) $\Rightarrow$ 3) $\Rightarrow$ 4) $\Rightarrow$ 5): It is obvious.
(Note that the assumption that $G$ is not elementary abelian is needed for the
implication 2) $\Rightarrow$~3).)

3) $\Rightarrow$ 1): This is nontrivial; it follows from [2, Theorem~1].

5) $\Rightarrow$ 4): Let $g\in G$ and $E\in \cE_G$.
For $x_E\in R$ such that $N_E(x_E) = 1$, define
$$x_{gEg^{-1}} = g(x_E). \eqno (1.3)$$
Then
$$N_{gEg^{-1}}(x_{gEg^{-1}}) = N_{gEg^{-1}}(g(x_E))
= g\bigl( N_E(x_E) \bigr) = g(1) = 1.$$

4) $\Rightarrow$ 3): Any $E\in \cE_G$ is a subgroup of an element
of~$\cE_G^{\rm max}$.
Then proceed as for 1) $\Rightarrow$ 2).
\hfill\cqfd
\medskip

It can be seen from (1.2) and (1.3) that the number of variables in
a formula $\Phi_G$ for $G$ can be reduced; we can restrict ourselves
to the variables $g(x_E)$, where $g\in G$ and where $E \in \cE_G^{\rm max}$
or $E \in \cE_G^0$.

\bigskip\goodbreak
\noindent
{\bf 2.~Three reductions}
\medskip

In this section we reduce in three steps the problem of finding a formula
for~$G$
to the problem of finding formulas for smaller groups of a special type.
\medskip

\noindent
{\it First reduction}
\smallskip

We start by reducing the problem to $p$-groups, where $p$ is a prime number.
Given a group~$G$, let $n = p_1^{a_1} \cdots p_r^{a_r}$
be the factorization of the order~$n$ of~$G$ in prime factors,
where $p_1, \ldots, p_r$ are distinct prime numbers, $r\geq 2$,
and the exponents $a_1, \ldots, a_r$ are positive integers.
Choose integers $d_1, \ldots, d_r$ such that
$$d_1 \, {n^{\;} \over p_1^{a_1} }+ \cdots + d_r\,  {{n^{\;} \over p_r^{a_r}}} = 1.$$
For every $i = 1, \ldots, r$, let $S_i$ be a Sylow $p_i$-subgroup
(of order~$p_i^{a_i}$) of~$G$.

The following result implies that, in order to find a formula for a group $G$,
it~is sufficient to find a formula for a Sylow $p$-subgroup of~$G$
for each prime number~$p$ dividing the order of~$G$.

\medskip\goodbreak
\noindent
{\sc Proposition~2.1}.---
{\it For each $i = 1, \ldots, r$, let $\Phi_{S_i}$ be a formula for~$S_i$. If
$$\Phi_G = d_1 \Phi_{S_1} + \cdots + d_r \Phi_{S_r},$$
then $\Phi_G$ is a formula for~$G$.
}
\medskip

\Pr
Suppose we are given a ring $R$ on which $G$ acts and
elements $x_E\in R$ such that $N_E(x_E) = 1$, one for each
elementary abelian subgroup $E$ of~$G$.
Replacing each variable $x_E$ ($E\in \cE_{S_i} \subset \cE_G$)
in the polynomial~$\Phi_{S_i}$ by the element $x_E\in R$,
we obtain an element $x_i \in R$
such that $N_{S_i}(x_i) = 1$ for each $i = 1, \ldots, r$.
Let us check that $N_G(x_G) = 1$ for $x_G = d_1 x_1 + \cdots + d_r x_r$.
Indeed,
$$N_G(x_i) = \sum_{g\in G}\, g(x_i) =
\sum_{g\in G/S_i}\, g \bigl( N_{S_i}(x_i) \bigr)
= \sum_{g\in G/S_i}\, g(1) = {n/p_i^{a_i}}.$$
Consequently,
$N_G(x_G) = d_1 {n/p_1^{a_1}} + \cdots + d_r {n/p_r^{a_r}} = 1$.
\hfill\cqfd
\goodbreak\medskip

Let us illustrate Proposition~2.1 in the case
of the symmetric group~$S_3$.
Let $s$ be a transposition and $t$ be a cyclic permutation in~$S_3$.
The subgroups $E_s$ and $E_t$ generated respectively by $s$ and~$t$
are elementary abelian. Then by Proposition~2.1,
$$\Phi_{S_3} = x_{E_s} - x_{E_t}$$
is a formula for~$S_3$.
\medskip

\goodbreak\noindent
{\it Second reduction}
\smallskip

We next reduce the problem from arbitrary $p$-groups to $p$-groups
that are extraspecial or almost extraspecial. Recall that
a $p$-group $G$ is {\it extraspecial}
(respectively {\it almost extraspecial}) if
$G$ fits into a central extension of the type
$$1 \to C_p \to G \to C_p^r \to 1,$$
where $r\geq 1$, and the center of $G$ is isomorphic to~$C_p$
(respectively to~$C_{p^2}$).
For a complete description of (almost) extraspecial groups,
see for instance~[5] (see also [8, Chapter~5]).
Note that with this definition no abelian group is extraspecial
and that the only abelian almost extraspecial group is~$C_{p^2}$.

For any $p$-group $G$ we denote $\cF_G$ the set of isomorphism classes of
groups $U$ satisfying the following conditions :

(i) $U$ is a subquotient (i.e., a homomorphic image of a subgroup) of~$G$ and

(ii) $U$ is extraspecial or almost extraspecial.

\noindent
We have $\cF_H \subset \cF_G$ whenever $H$ is a subquotient of~$G$.

The following result states that, in order to find a formula for a finite
$p$-group~$G$,
it is sufficient to have formulas for all groups in~$\cF_G$.

\medskip\goodbreak
\noindent
{\sc Theorem~2.2}.---
{\it For any finite $p$-group $G$ there is an algorithm whose output is a
formula $\Phi_G$ for~$G$ and whose inputs are formulas $\Phi_H$
for all $H \in \cF_G$.
}
\medskip

Before we prove the theorem, we establish two intermediate results.

\medskip\goodbreak
\noindent
{\sc Lemma~2.3}.---
{\it If a $p$-group $G$ is neither elementary abelian, nor extraspecial,
nor almost extraspecial,
then there is a central element $h\in G$ of order $p$ such that the quotient
group $G/\langle h\rangle$ is not an elementary abelian group.
}
\medskip

\Pr
If $G$ is abelian, then $G$ is of order $\geq p^3$.
Take an element $g$ of order $p$ in~$G$.
If $G/\langle g\rangle$ is not elementary abelian, we are done.
If $G/\langle g\rangle$ is elementary abelian, say $G/\langle g\rangle
\cong C_p^r$ for
some $r$ that is necessarily at least~$2$, then
by the classification of finite abelian $p$-groups we have
$G \cong C_{p^2} \times C_p^{r-1}$. Let $\sigma$ be an element of order $p$
in~$C_p^{r-1}$ (it exists since $r-1 \geq 1$) and $h\in G$
be the element mapped to $(0,\sigma) \in C_{p^2} \times C_p^{r-1}$.
Then $G/\langle h\rangle$ contains an element of order~$p^2$, hence
is not elementary abelian.

Now assume that $G$ is not abelian. Let $g$ be an element of order $p$
in the center $Z(G)$ of~$G$. If $G/\langle g\rangle$ is not elementary
abelian, we
are done. Therefore we may assume that $G/\langle g\rangle$ is elementary
abelian.
Observe that under this condition the commutator subgroup $G'$ of $G$
is the subgroup $\langle g\rangle$ generated by~$g$.
Since $G$ is neither abelian, nor extraspecial, nor almost extraspecial,
its center $Z(G)$ is not isomorphic to $C_p$ or to~$C_{p^2}$.
Moreover, $Z(G)$ is not cyclic of order $\geq p^3$ since
$G/\langle g\rangle$ is elementary abelian.
Therefore there is an element $h\in Z(G)$ of order $p$ such that
$\langle h\rangle \neq \langle g\rangle = G'$, and so $G/\langle h\rangle$
is not abelian, hence not elementary abelian.
\line{\hfill\cqfd}
\medskip

Let $G$ be a $p$-group that is neither elementary abelian, nor extraspecial,
nor almost extraspecial. By Lemma~2.3 there is a subgroup~$U$ of order~$p$
in the center of~$G$ such that $G/U$ is not elementary abelian.
We fix such a subgroup~$U$.

Let $\pi : G\to G/U$ be the natural projection.
Let $\Phi_{G/U}$ be a formula for $G/U$, and for each $\bar{E} \in \cE_{G/U}$
let $\Phi_{\pi^{-1}(\bar{E})}$ be a formula for the proper subgroup
$\pi^{-1}(\bar{E})$ of~$G$
(it is a proper subgroup because $G/U$ is not elementary abelian).
Set
$$\Phi_G = \Phi_{G/U}\bigl( N_U(\Phi_{\pi^{-1}(\bar{E})}) \bigr) x_U .
\eqno (2.1)$$
Equality~(2.1) defines a noncommutative polynomial with integer coefficients
in the variables $g(x_E)$, where $g\in G$ and $E\in \cE_G$.
This is a consequence of the following observations
on the right-hand side of~(2.1).

Firstly, $N_U(\Phi_{\pi^{-1}(\bar{E})})$ has the following meaning:
we replace each monomial
$h_1(x_{E_1})  \cdots h_s(x_{E_s})$
in~$\Phi_{\pi^{-1}(\bar{E})}$,
where $h_1, \ldots, h_s \in \pi^{-1}(\bar{E})$ and
$E_1, \ldots, E_s \in \cE_{\pi^{-1}(\bar{E})}$,
by the polynomial
$$\sum_{u\in U}\,  (uh_1)(x_{E_1})  \cdots (uh_s)(x_{E_s}).$$

Secondly, the expression $\Phi_{G/U}\bigl( N_U(\Phi_{\pi^{-1}(\bar{E})})
\bigr)$
means that we replace each letter $x_{\bar{E}}$ ($\bar{E} \in \cE_{G/U}$)
in the polynomial~$\Phi_{G/U}$ by the polynomial
$N_U(\Phi_{\pi^{-1}(\bar{E})})$ whose meaning has just been explained.
In this way, each variable $\bar{g}(x_{\bar{E}})$, where $\bar{g} \in G/U$
and $\bar{E} \in \cE_{G/U}$, becomes a polynomial in the variables~$g(x_E)$,
where $g\in G$ and
$E \in \cup_{\bar{E} \in \cE_{G/U}}\, \cE_{\pi^{-1}(\bar{E})} \, (\subset
\cE_G)$.

Finally, the polynomial
$\Phi_{G/U}\bigl( N_U(\Phi_{\pi^{-1}(\bar{E})}) \bigr)$ is multiplied
on the right by the variable~$x_U$, which makes sense since $U\in \cE_G$.

\medskip\goodbreak
\noindent
{\sc Proposition~2.4}.---
{\it With the previous notation, $\Phi_G$ is a formula for~$G$.
}
\medskip

\Pr
Suppose $G$ acts on a ring~$R$ and we have elements $x_E \in R$ such that
$N_E(x_E) = 1$, one for each~$E\in \cE_G$.
In particular, we have an element $x_U \in R$ such that $N_U(x_U) = 1$.

For each $\bar{E} \in \cE_{G/U}$, let $x_{\pi^{-1}(\bar{E})}$ be the
element of~$R$ obtained from the formula $\Phi_{\pi^{-1}(\bar{E})}$
by replacing each variable $h(x_E)$,
where $h\in \pi^{-1}(\bar{E})$ and $E \in \cE_{\pi^{-1}(\bar{E})}
\,(\subset \cE_G)$,
by the value of the action of~$h$ on the element~$x_E\in R$.
By definition of a formula, we have
$$N_{\pi^{-1}(\bar{E})}(x_{\pi^{-1}(\bar{E})}) = 1.$$
The element $N_U(x_{\pi^{-1}(\bar{E})})$ clearly belongs to the subring~$R^U$.
Let
$$z_{G/U} = \Phi_{G/U}\bigl( N_U(x_{\pi^{-1}(\bar{E})}) \bigr)$$
be obtained by replacing each variable $\bar{g}(x_{\bar{E}})$ of~$\Phi_{G/U}$
by $\bar{g}N_U(x_{\pi^{-1}(\bar{E})}) \in R^U$
(this makes sense since $G/U$ acts on~$R^U$).
Since $\Phi_{G/U}$ is a formula for~$G/U$, we have
$N_{G/U}(z_{G/U}) = 1$. The element $z_{G/U}$ belongs to~$R^U$ because
the inputs in its definition are in this subring.
To conclude, let $x_G = z_{G/U} x_U \in R$.
Using the $R^U$-linearity of~$N_U$, we obtain
$$\eqalign{
N_G(x_G) & = N_{G/U}(N_U(z_{G/U} x_U)) \cr
& = N_{G/U}(z_{G/U}N_U(x_U)) \cr
& = N_{G/U}(z_{G/U} \cdot 1) = 1. \cr
}$$
\hfill\cqfd
\medskip

\goodbreak\noindent
{\sc Proof of Theorem~2.2}.---
We proceed by induction on the order of~$G$.

Suppose $G$ is of order~$p^3$.
If $G$ is elementary abelian or extraspecial, we are done
(note that $G$ cannot be almost extraspecial).
Otherwise, $G \cong C_{p^3}$ or $G \cong C_{p^2} \times C_p$, and in
both cases $\cF_G = \{ C_{p^2}\}$.
Fix a central subgroup~$U$ of order~$p$
such that $G/U$ is not elementary abelian. Then $G/U$, being of order~$p^2$,
is isomorphic to~$C_{p^2}$, which is almost extraspecial. 
The group $G/U$ has a unique elementary abelian subgroup~$\bar{E}$ of order~$p$
whose lifting $\pi^{-1}(\bar{E})$ to~$G$, being of order~$p^2$,
is either isomorphic to~$C_p^2$ (elementary abelian)
or~$C_{p^2}$ (belonging to~$\cF_G$).
In both cases, by Proposition~2.4, (2.1) yields a formula for~$\Phi_G$
in terms of~$\Phi_{C_{p^2}}$.

Let $G$ be a group of order~$p^n$ with $n \geq 4$.
Suppose we have proved the theorem for all groups of order at most~$p^{n-1}$.
As above, we may assume that $G$ is neither elementary abelian,
nor extraspecial, nor almost extraspecial.
We again fix a central subgroup~$U$ of order~$p$
such that $G/U$ is not elementary abelian.
By the induction hypothesis, there is an algorithm whose output is a
formula $\Phi_{G/U}$ for~$G/U$
and whose inputs are formulas for all groups in~$\cF_{G/U}$.
The lifting $\pi^{-1}(\bar{E})$ to~$G$ of each $\bar{E} \in \cE_{G/U}$
is a proper subgroup of~$G$ since $G/U$ is not elementary abelian.
By the induction hypothesis again, there is an algorithm whose output is a
formula $\Phi_{\pi^{-1}(\bar{E})}$ for~$\pi^{-1}(\bar{E})$ and
whose inputs are formulas for all groups in~$\cF_{\pi^{-1}(\bar{E})}$.
Therefore by Proposition~2.4, (2.1) yields an algorithm whose output is a
formula $\Phi_G$ for~$G$ and whose inputs are formulas
for all groups in the sets~$\cF_{\pi^{-1}(\bar{E})}$ or in~$\cF_{G/U}$.
We conclude by observing that $\cF_{G/U} \subset \cF_G$
and $\cF_{\pi^{-1}(\bar{E})} \subset \cF_G$ for all~$\bar{E} \in \cE_{G/U}$.
\hfill\cqfd

\medskip
\goodbreak\noindent
{\it Third reduction}
\smallskip

We now consider the case when $G$ is the product of two groups.

\medskip
\noindent
{\sc Theorem~2.5}.---
{\it Let $G = G_1 \times G_2$ be a product of two $p$-groups $G_1$ and~$G_2$.
There is an algorithm whose output is a formula $\Phi_G$ for~$G$
and whose inputs are formulas for all groups in~$\cF_{G_1} \cup\cF_{G_2}$.
}
\medskip

Before we give the proof of this theorem, we establish a result
similar to Proposition~2.4.

Consider the case when $G = H \times C$, where $H$ is a $p$-group and
$C$ is an elementary abelian $p$-group of rank one, i.~e.,~$C \cong C_p$.
Set
$$\Phi_G = \Phi_H\bigl( N_C(x_{E'\times C})\bigr) x_C, \eqno (2.2)$$
where $\Phi_H$ is a formula for $H$ and $E'\in \cE_H \, (\subset \cE_G)$.
Observe that $E' \times C \in \cE_G$. The right-hand side of~(2.2)
has the following meaning. Firstly,
$$N_C(x_{E'\times C}) = \sum_{u\in C}\, u(x_{E'\times C}).$$
Secondly, $\Phi_H\bigl( N_C(x_{E'\times C})\bigr)$ means that we
replace each letter $x_{E'}$ ($E' \in \cE_H$)
in the polynomial~$\Phi_H$ by the polynomial
$N_C(x_{E'\times C})$ defined above.
In this way, each variable $h(x_{E'})$ of~$\Phi_H$, where $h \in H$
and $E' \in \cE_H$, becomes a polynomial in the variables~$g(x_{E'\times C})$,
where $g\in G$ and $E'\times C \in \cE_G$.
Therefore, the right-hand side of~(2.2) is a noncommutative polynomial
with integer coefficients and in the right set of variables.

\medskip\goodbreak
\noindent
{\sc Proposition~2.6}.---
{\it With the previous notation, $\Phi_G$ is a formula for~$G$.
}
\medskip

\Pr
Suppose $G$ acts on a ring $R$ and we have elements $x_E \in R$ such that
$N_E(x_E) = 1$, one for each $E\in G$.
In particular, we have $x_C \in R$ such that $N_C(x_C) = 1$.
For each $E' \in \cE_H$, the product $E'\times C$ is an
elementary abelian subgroup of~$G$, and we have an element
$x_{E'\times C} \in R$ such that $N_{E'\times C}(x_{E'\times C}) = 1$.
The element $N_C(x_{E'\times C})$ belongs to $R^C$ and satisfies
$$N_{E'}\bigl( N_C(x_{E'\times C}) \bigr) = 1.$$
Let
$z = \Phi_H\bigl( N_C(x_{E'\times C})\bigr)$
be the element obtained by replacing each variable $h(x_{E'})$ of $\Phi_H$
by $hN_C(x_{E'\times C}) \in R^C$ (the group $H = G/C$ acts on~$R^C$).
Since $\Phi_H$ is a formula for~$H$, we have $N_H(z) = 1$.
Moreover, since the inputs belong to~$R^C$, so does~$z$.
Now let $x_G = z x_C$.
Then, using the $R^C$-linearity of~$N_C$, we obtain
$$\eqalign{
N_G(x_G) & = N_H(N_C(z x_C)) \cr
& = N_H(z N_C(x_C)) \cr
& = N_H(z \cdot 1) = 1. \cr
}$$
\hfill\cqfd
\medskip

\goodbreak
\noindent
{\sc Proof of Theorem~2.5}.---
(a) Assume first that $G_2$ is elementary abelian.
Let us prove by induction on the order of~$G_2$ that there is
an algorithm whose output is a formula $\Phi_G$ for~$G$
and whose inputs are formulas for all groups in~$\cF_{G_1}$.
Since $\cF_{G_2} = \emptyset$, it will prove Theorem~2.5 in this case.

If $G_2$ is of order~$p$, we appeal to Proposition~2.6:
Formula~(2.2), in which we have replaced $H$ by $G_1$ and $C$ by $G_2$,
yields an algorithm whose output is a formula $\Phi_G$ for the group~$G$
and whose input is a formula for~$G_1$. Therefore, by Theorem~2.2
there is an algorithm whose output is a formula for~$G$ and whose
inputs are formulas for all groups in~$\cF_{G_1}$.

If $G_2$ is of order $> p$, we write $G_2 = G'_2 \times C$, where $C$ is
cyclic of
order~$p$ (the subgroup~$G'_2$ is elementary abelian).
Reasoning as above, we obtain an algorithm
whose output is a formula for~$G = G_1 \times G'_2 \times C$ and whose
input is a formula for~$G_1 \times G'_2$.
By induction there is an algorithm whose output is a formula for~$G_1
\times G'_2$
and whose inputs are formulas for all groups in~$\cF_{G_1}$.
Therefore there is an algorithm whose output is a formula for~$G$
and whose inputs are formulas for all groups in~$\cF_{G_1}$.

(b) Let $G_2$ be an arbitrary $p$-group of order~$p^n$ with $n\geq 2$
and suppose we have proved Theorem~2.5 for all groups $G = G_1 \times G'_2$
such that the order of~$G'_2$ is at most~$p^{n-1}$.
By Part~(a) we may assume that $G_2$ is not elementary abelian.

Let $U$ be a central subgroup of~$G_2$ of order~$p$ and
$$\pi : G = G_1 \times G_2 \to G/U = G_1 \times G_2/U$$
be the natural projection.
By Proposition~2.4, Formula~(2.1) yields an algorithm
whose output is a formula for~$G$
and whose inputs are formulas for~$G/U$ and for $\pi^{-1}(\bar{E})$,
where $\bar{E}$ runs over~$\cE_{G/U}$.
By induction there is an algorithm whose output is a formula
for~$G/U = G_1 \times G_2/U$ and whose inputs are formulas
for all groups in~$\cF_{G_1} \cup \cF_{G_2/U}$.

Each elementary abelian subgroup $\bar{E}$ of $G/U = G_1 \times G_2/U$
is clearly contained in an  elementary abelian subgroup
of the form~$E_1 \times E_2$,
where  $E_1$ is an elementary abelian subgroup of~$G_1$ and
$E_2$ is an elementary abelian subgroup of~$G_2/U$.
Its inverse image $\pi^{-1}(\bar{E})$ is therefore contained
in a group of the form $E_1 \times N$, where $N$ is a subgroup of~$G_2$.
Formula~(1.2) shows how to obtain a formula for~$\pi^{-1}(\bar{E})$
from a formula for~$E_1 \times N$.
By Part~(a) there is an algorithm whose output is a formula for~$E_1 \times N$
and whose inputs are formulas for all groups in~$\cF_N$.
Summing up and observing that $\cF_{G_2/U} \subset \cF_{G_2}$ and
$\cF_N \subset \cF_{G_2}$,
we conclude that there is algorithm whose output is a formula for~$G$
and whose inputs are formulas for all groups in~$\cF_{G_1} \cup \cF_{G_2}$.
\hfill \cqfd
\medskip

\medskip
\noindent
{\sc Remarks~2.7}.
(a) Lemma~2.3 can be used to show that $\cF_G = \emptyset$
if and only if $G$ is elementary abelian.

(b) By definition of~$\cF_G$, if $G$ is an {\it abelian} $p$-group,
then $\cF_G$ is empty (if $G$ is elementary abelian) or $\cF_G = \{ C_{p^2}
\}$.
By Theorem~2.2 it suffices to have a formula for $C_{p^2}$ in order to
obtain formulas for all abelian groups. A formula for~$C_{p^2}$
was obtained in [3, Corollary~1].
We recall it here for the sake of completeness:
let $\sigma$ be a generator of~$C_{p^2}$ and $E$ be the unique abelian
elementary subgroup of~$C_{p^2}$. Then
$$\eqalign{
\Phi_{C_{p^2}} & = x_E^2 +
\sum_{j = 0}^{p-1}\sum_{k=1}^{p -1} \sum_{i=0}^{k -1} \,
\sigma^{ip}( x_E) \,\sigma^{j-(k-i)p}(x_E) \, x_E \cr
& {} \hskip 12pt - \sum_{j = 0}^{p-1}\sum_{k=1}^{p -1} \sum_{i=0}^{k -1}\,
\sigma^{ip+1}( x_E) \,\sigma^{j-(k-i)p+1}(x_E) \, x_E \cr
& {} \hskip 12pt - \sum_{k=1}^{p -1} \sum_{i=0}^{k -1}\,
\sigma^{ip}( x_E)  \, x_E
+ \sum_{k=1}^{p -1} \sum_{i=0}^{k -1}\,
\sigma^{ip+1}( x_E)  \, x_E \cr
} \eqno (2.3)$$
is a formula for $C_{p^2}$. Formula~(0.2) is the special case of~(2.3)
when~$p=3$.

(c) Theorem~2.5 may help get a better reduction than Theorem~2.2 for
nonabelian groups.
For instance, take the product $G= Q_8\times Q_8$ of two copies of
the quaternion group of order~$8$.
By~[8, Chapter~5] the central product $G_1$ of $Q_8$ with itself,
which is a quotient of~$G$, is an extraspecial group;
it belongs to~$\cF_G$.
Theorem~2.2 therefore suggests that a formula for the group $G_1$  (of
order~32) is needed
to obtain a formula for~$G$.
Nevertheless, by Theorem~2.5 only formulas for the groups
in $\cF_{Q_8} = \{ C_4, Q_8\}$ are needed.

\bigskip\goodbreak
\noindent
{\bf 3.~Cohomological preliminaries}
\medskip

In this section we present two results needed in the sequel.
The first one is an important consequence of the existence of
a norm one element.

\medskip
\noindent
{\sc Proposition~3.1}.---
{\it Let $G$ be a finite $p$-group acting on a ring~$R$. If there is
an element $x\in R$ such that $N_G(x) = 1$, then
$H^i(G,R) = 0$ for all~$i>0$.
}
\medskip

\Pr
We proceed by induction on the order of~$G$.
Assume first that $G$ is a cyclic group of order~$p$ with generator~$\sigma$.
Recall that the cohomology groups of a cyclic group are given for all
$j\geq 1$ by
$$H^{2j}(G,R) = R^G/\Im\ N_G \and H^{2j-1}(G,R) = \Ker\ N_G/(\sigma - 1)(R).$$
The even cohomology groups vanish
since $N_G(R^G x) = R^GN_G(x) = R^G$.
The vanishing of the odd cohomology groups follows from [3, Lemma~1].

Now let $G$ be of order $p^n$ with $n\geq 2$ and assume that the lemma holds
for all $p$-groups of order~$< p^n$. Take a normal subgroup $U$ of~$G$
of index~$p$ (such a subgroup always exists).
Using Formula~(1.2), out of the element $x\in R$ satisfying $N_G(x) = 1$
we derive an element $x_U \in R$ such that $N_U(x_U) = 1$.
The existence of~$x_U$ implies by induction that $H^i(U,R) = 0$ for all~$i>0$.

It then follows from the Lyndon-Hochschild-Serre spectral sequence
that the inflation maps
$$\Inf : H^i(G/U, R^U) \to H^i(G,R)$$
are isomorphisms for all $i$. Now $G/U$ is cyclic and the element
$N_U(x)\in R^U$
satisfies
$$N_{G/U}(N_U(x)) = N_G(x) = 1.$$
Then by the first part of the proof the cohomology groups $H^i(G/U, R^U)$
vanish
for all~$i > 0$, and so do the groups~$H^i(G,R)$.
\hfill\cqfd
\medskip

In the next sections we shall represent elements of the cohomology
group $H^1(G,M)$, where $G$ is a group and $M$ is a left $G$-module,
by 1-cocycles of~$G$ with values in~$M$.
Recall from [6, Chapter~X, \S~4]
that such a 1-cocycle (also called a crossed homomorphism)
is a map $\beta : G \to M$ satisfying
$$\beta(gh) = \beta(g) + g\beta(h) \eqno (3.1)$$
for all $g$, $h\in G$.
A~1-coboundary is a map $\beta : G \to M$ for which there
exists $m\in M$ such that
$$\beta(g) = (g - 1)m  \eqno (3.2)$$
for all $\in G$. It is easy to check that a 1-coboundary is a 1-cocycle.
1-cocycles and 1-coboundaries are elements of the standard cochain complex
whose cohomology is~$H^*(G,M)$.
If $\delta$ is the differential in the standard cochain complex,
then (3.2) can be rewritten as~$\beta = \delta(m)$.

The following identities are easy consequences of the functional
equation~(3.1).

\medskip\goodbreak
\noindent
{\sc Lemma~3.2}.---
{\it Let $\beta : G \to M$ be a 1-cocycle of~$G$ with values in
a left $G$-module~$M$.

(a) For the neutral element $1\in G$, we have $\beta(1) = 0$.

(b) For $g\in G$ and $i \geq 2$, we have
$$\beta(g^i) = (1 + g + \cdots + g^{i-1}) \beta(g).$$

(c) If $g\in G$ is of order~$N$, then
$$(1 + g \cdots + g^{N-1})\beta(g) = 0.$$

(d) For any $g\in G$, we have $\beta(g^{-1}) = - g^{-1} \beta(g)$.

(e) If $\sigma$, $\tau \in G$ satisfy $\tau\sigma = \sigma^{-1}\tau$, then
$$(\sigma - 1) \beta(\tau) + (1 + \sigma\tau) \beta(\sigma) = 0.$$
}

\bigskip\goodbreak
\noindent
{\bf 4.~A method for finding formulas for $p$-groups}
\medskip

We now present a method for finding a formula for
a given $p$-group~$G$. It consists in taking a presentation of~$G$
and deriving from it a system of equations
whose indeterminates are elements $b(\sigma) \in R$, one for each generator
$\sigma$
in the presentation. There is an equation for each relation in
the presentation.
Group cohomology guarantees that this system of equations has a solution.
Once we have a solution, we again use homological algebra
to obtain an explicit formula for~$G$.

In view of Propositions 1.3 and~2.1 we may assume that
$G$ is a finite $p$-group (not elementary abelian) and that we have
formulas for all proper subgroups of~$G$.
Let $G$ act on a ring~$R$.
We assume the existence of $x_H\in R$ such that $N_H(x_H) = 1$ for each
proper subgroup~$H$ of~$G$.
Our aim is to give an explicit formula for $x_G \in R$ with $N_G(x_G) = 1$
in terms of the elements $x_H$ and of the elements of~$G$.

Fix a normal subgroup $U$ of index~$p$ in~$G$ and
choose an element $\sigma\in G$
whose class $\bar{\sigma}$ generates the cyclic group~$G/U$.
Set $x = x_U \in R$ (this is one of the elements $x_H$ whose
existence was assumed above); we have~$N_U(x) = 1$.

\medskip\goodbreak
\noindent
{\sc Proposition~4.1}.---
{\it Let $a\in R$ be a $U$-invariant element such that
$$N_{G/U}(a) = (1 + \sigma + \cdots + \sigma^{p-1})(a) = 1.$$
Then $N_G(y) = 1$ if $y = ax$ or $y = xa$.
}
\medskip

\Pr
Let $y = ax$. By the $R^U$-linearity of~$N_U$, we have
$$N_G(y) = N_{G/U}(N_U(ax)) = N_{G/U}(a N_U(x)) = N_{G/U}(a) = 1.$$
A similar proof holds for $y = xa$.
\hfill\cqfd
\medskip

To solve the problem for~$G$, it is therefore sufficient to
find an element $a\in R^U$ such that $N_{G/U}(a) = 1$.
In the rest of the section we show how to find such an element.

We start as in~[3, Section~2] by considering
the group $B = \Hom_{\ZZ}(\ZZ[G],R)$ of $\ZZ$-linear maps from the
group ring~$\ZZ[G]$ to~$R$. The group $G$ acts on the left on~$B$
by $(g\varphi)(s)  = \varphi(sg)$
for all $g,s\in G$ and $\varphi \in B$.
The ring $R$ is a $G$-submodule of~$B$, where an element $r\in R$ is
identified with
the element $\varphi_r \in B$ given by $\varphi_r(g) = g(r)$
for all $g\in G$. Let $C = B/R$ be the quotient $G$-module. We denote
$q : B\to C$ the canonical surjection.

Consider the following commutative square:
$$\matrix{
H^1(G/U,C^U) & \hfl{\Inf}{} & H^1(G,C)\cr
\noalign{\smallskip}
\vfl{\delta}{} & & \vfl{\delta}{} \cr
\noalign{\smallskip}
H^2(G/U,R^U) & \hfl{\Inf}{} & H^2(G,R)\cr
} \eqno (4.1)
$$
The vertical maps~$\delta$ in (4.1) are the connecting homomorphisms
arising from the short exact sequences
$$ 0 \to R \to B \to C \to 0 \and
0 \to R^U \to B^U \to C^U \to H^1(U,R) = 0.$$
By Proposition~3.1 the group $H^1(U,R)$ vanishes because of the existence
of the element~$x \in R$ satisfying~$N_U(x) = 1$.
The maps $\delta$ are isomorphisms because
$B$ is a co-induced $G$-module and $B^U$ is a co-induced $G/U$-module,
hence $H^i(G,B) = H^i(G/U, B^U) = 0$ for all~$i > 0$.

The horizontal maps in (4.1) are inflation maps.
By the Lyndon-Hochschild-Serre spectral sequence
the vanishing of $H^i(U,R)$ for~$i>0$ (see Proposition~3.1)
implies that the lower inflation map
$\Inf : H^2(G/U,R^U) \to H^2(G,R)$
is an isomorphism. Therefore, the upper inflation map is an isomorphism as
well.
In view of this, of Proposition~1.3 and of Proposition~3.1,
all groups in the square~(4.1) vanish.

Next, define $\varphi \in B$ by
$$\varphi(g) = \cases{ 1 & if $g \in U$,
\cr
\noalign{\smallskip}
$0$ & otherwise.\cr}
\eqno (4.2)$$
It is clear that $\varphi$ is invariant
under the action of the subgroup~$U$ and that
$$N_{G/U}(\varphi) =
(1 + \sigma + \cdots + \sigma^{p-1})(\varphi) = \varphi_1 = 1.$$
Consider the map $\alpha_0 : G/U\to B^U$ given by
$$\alpha_0(\bar{\sigma}^k) = \cases{
$0$ & if $\, k=0$, \cr
\noalign{\smallskip}
\varphi & if $\, k=1$, \cr
\noalign{\smallskip}
(1 + \sigma + \cdots + \sigma^{k-1})(\varphi) & if $\, 2 \leq k \leq p-1$.\cr
}
\eqno (4.3)$$
Let $q\alpha_0 : G/U \to C^U$ be the composition of $\alpha_0$
with $q_{| B^U}: B^U\to C^U$.

\medskip\goodbreak
\noindent
{\sc Lemma~4.2}.---
{\it The map $q\alpha_0 : G/U \to C^U$ is a 1-cocycle.
}
\medskip

\Pr
It suffices to check that
$$q\alpha_0(\bar{\sigma}^i) + \bar{\sigma}^i q\alpha_0(\bar{\sigma}^j) =
q\alpha_0(\bar{\sigma}^{i+j})$$
for all $i, j\in \{0, 1, \ldots, p-1\}$.
This follows from the definition of~$\alpha_0$ and the
following equalities in~$C$:
$$q((1 + \bar{\sigma} + \cdots + \bar{\sigma}^{p-1})(\varphi))
 = q(N_{G/U}(\varphi)) = q(\varphi_1) = 0. \eqno \hbox{\cqfd}$$

Define $\alpha: G\to B$ by $\alpha(g) = \alpha_0(\bar{g})$
for all $g\in G$, where $\bar{g}$ denotes the class of~$g$ in~$G/U$.
By~(4.3) the value of $\alpha$ on the chosen element $\sigma$
is $\alpha(\sigma) = \varphi$.

Let $q\alpha : G \to C$ be the composition of $\alpha$ with $q: B\to C$.
By Lemma~4.2, $q\alpha_0 : G/U \to C^U$ represents an element
in~$H^1(G/U,C^U)$.
It is easy to check that its image in $H^1(G,C)$ under the upper inflation map
in the square~(4.1) is represented by the map~$q\alpha : G \to C$.
Since $q\alpha_0$ is a 1-cocycle, so is~$q\alpha$.

Our idea is to correct $\alpha : G \to B$ as follows.

\medskip\goodbreak
\noindent
{\sc Lemma~4.3}.---
{\it There is a map $b : G\to R$ such that
$$\alpha - b = \alpha - \varphi_b : G\to B$$
is a $1$-cocycle with values in~$B$.
}
\medskip

\Pr
We have seen above that $H^1(G,C) = 0$.
Since $q\alpha$ is a $1$-cocycle with values in~$C$, it is a $1$-coboundary;
so there is $\bar\psi \in C$ such that $q\alpha = \delta(\bar\psi)$.
Lift $\bar\psi$ to an element $\psi \in B$ and set $b = \alpha - \delta(\psi)$.
Then $\alpha - b = \delta(\psi)$ is a 1-coboundary with values in~$B$,
hence a 1-cocycle.
\hfill\cqfd
\medskip

We now claim that it suffices to perform the following three tasks
in order to find a formula for~$G$.

\medskip\goodbreak
\noindent
{\bf Task~1:} Write the set of equations satisfied by the values $b(g)\in R$
of the map~$b$, obtained by expressing that $\alpha - b : G\to B$ is a
$1$-cocycle.
We can reduce the number of unknowns by choosing a presentation
$\langle \sigma_1, \ldots, \sigma_r \, |\, R_1, \ldots, R_s\rangle$ of~$G$.
Since the generators are of finite order, we may assume that each relation
$R_i$ is a word in $\sigma_1, \ldots, \sigma_r$ (i.e., the inverses
of $\sigma_1, \ldots, \sigma_r$ do not appear in~$R_i$).

Set $\beta = \alpha - b$.
For each relation $R_i = \sigma_{i_1} \sigma_{i_2} \cdots \sigma_{i_t}$
define
$$\beta(R_i)
= \beta(\sigma_{i_1}) + \sigma_{i_1} \beta(\sigma_{i_2}) + \cdots
+ \sigma_{i_1} \sigma_{i_2}  \cdots \sigma_{i_{t-1}} \beta(\sigma_{i_t}).
\eqno (4.4)$$
By setting $\beta(R_i) = 0$ for all $i = 1, \ldots, s$, we obtain
a system $(\Sigma)$ of $s$ equations whose unknowns are
$b(\sigma_1), \ldots, b(\sigma_r)$.
It is an easy exercise to show that the values $b(\sigma_1), \ldots,
b(\sigma_r)\in R$
determine uniquely a map $b: G\to R$ such that $\alpha - b$ is a 1-cocycle.

\medskip\goodbreak
\noindent
{\bf Task~2:} By Lemma~4.3 the system of equations $(\Sigma)$
derived in Task~1 has a solution $b: G\to R$.
Task~2 consists in writing down such a solution 
polynomially in terms of the given data.

\medskip\goodbreak
\noindent
{\bf Task~3:} By Proposition~3.1 the existence of a norm one element $x\in R$
for~$U$ implies the vanishing of~$H^1(U,R)$. Hence,
for any 1-cocycle $\beta : U\to R$ there is an element
$w\in R$ such that $\beta(g) = (g-1)w$ for all $g\in U$.
Give an explicit expression of such an element $w$ as a noncommutative
polynomial with integer coefficients 
in the variables $u(x)$ and $u(\beta(v))$, where $u,v \in U$.

\medskip
Once Tasks 1--3 are completed, we solve the problem for~$G$ as follows.
Let $b: G\to R$ be a solution of the system~$(\Sigma)$
(in particular, we have an element~$b(\sigma) \in R$).
Then $\alpha - b$ is a 1-cocycle with values in~$B$.
Since $B$ is cohomologically trivial,
there is $\psi : G\to R$ such that $\alpha - b = \delta(\psi)$.
In particular, since $\alpha(g) = \varphi$ for $g = \sigma$,
we obtain
$$\varphi - b(\sigma) = (\sigma - 1)\psi . \eqno (4.5)$$
Similarly, by (4.3),
$$0 - b(g) = (g - 1) \psi \eqno (4.6)$$
for all $g\in U$.
Equations~(4.6) imply that the restriction of~$b$ to~$U$ is a $1$-cocycle
with values in~$R$.
After performing Task~3, we have an element $w\in R$ such that
$$b(g)  = (g - 1) w \eqno (4.7)$$
for all $g\in U$.
Relations~(4.6--4.7) together imply
$$(g - 1)(\psi + w) = 0$$
for all $g\in U$, which means that $\psi + w$ is $U$-invariant.

\medskip\goodbreak
\noindent
{\sc Proposition~4.4}.---
{\it With the previous notation  the element
$$a  = b(\sigma) + (1 - \sigma)(w) \in R$$
is $U$-invariant and we have $N_{G/U}(a) = 1$.
}
\medskip

\Pr
(a) Relation~(4.5) allows us to express $a$ under the form
$$a = \varphi - (\sigma - 1)(\psi + w).$$
To check the $U$-invariance of~$a$,
it is enough to check the $U$-invariance of
$$(\sigma - 1)(\psi + w)$$
since $\varphi$ is $U$-invariant.
For $u\in U$ let $u'\in U$ be such that $u\sigma = \sigma u'$
(recall that $U$ is a normal subgroup of~$G$).
Therefore, in view of the $U$-invariance of $\psi + w$,
for $u\in U$ we obtain
$$\eqalign{
u(\sigma - 1)(\psi + w)
& = u\sigma (\psi+w) - u(\psi + w)\cr
& = \sigma u'(\psi+w) - u(\psi + w)\cr
& = \sigma (\psi + w) - (\psi + w) \cr
& = (\sigma - 1)(\psi + w). \cr
}$$

(b) Since $\sigma^p$ belongs to~$U$, the $U$-invariance of $\psi + w$ implies
$$\eqalign{
N_{G/U}(a) & = N_{G/U}(\varphi)
- (1 + {\sigma} + \cdots + {\sigma}^{p-1})(\sigma -1)(\psi+ w) \cr
& = 1 - (\sigma^p -1)(\psi + w) = 1. \cr }$$
\line{\hfill\cqfd}
\medskip

By Proposition~4.1, the element $y  = ax \in R$ (or $y = xa$) then
yields an explicit norm one element for~$G$ with the appropriate form.
This solves the problem for~$G$.

\medskip\goodbreak
\noindent
{\sc Remark~4.5}.
Before we close this section, let us evaluate the level of difficulty of Tasks~1--3.
We explained above how to perform Task~1; this is easy.
For Task~3 we have a general method to solve it; it will be detailed
in the next section.

At the moment we do not have a general method to solve Task~2,
which consists in solving the system of equations~$(\Sigma)$ defined above.
The solutions given in Sections~6--8 have been found in an {\it ad hoc} way;
we nevertheless observe that they are of a very simple form. 
If we could prove in full generality that the solutions of~$(\Sigma)$ are of the form
$$\sum_{H}\, A_H (x_H),$$
where $H$ runs over all maximal proper subgroups of~$G$, $A_H \in \ZZ[G]$,
and $x_H\in R$ satisfies $N_H(x_H) = 1$, then solving $(\Sigma)$
could be reduced to solving a system~$(\Sigma')$ of linear equations
over the ring of integers~$\ZZ$, whose number of unknowns and of equations
can be bounded in terms of~$G$. 
More precisely, if $r$ is the number of generators and $s$ is the number of relations
in the chosen presentation of the group~$G$,
and if $m$ is the number of maximal proper subgroups of~$G$, then
the number of unknowns in~$(\Sigma')$ is $rm|G|$ and
the number of equations is~$s |G|$. 
Note that the number of maximal proper subgroups of~$G$ is
$m = 
(p^N - 1)(p-1)$, where $p^N$ is the order of the quotient of~$G$
by its Frattini subgroup.

\bigskip\goodbreak
\noindent
{\bf 5.~Writing a 1-cocycle as an explicit 1-coboundary}
\medskip

We consider a finite $p$-group $U$ acting on a ring~$R$.
Assume that we have an element $x\in R$ such that $N_U(x) = 1$.
The cohomology group $H^1(U,R)$ vanishes by Proposition~3.1.
Therefore, given a 1-cocycle $\beta : U\to R$, there exists
$w\in R$ such that $\beta(g) = (g-1)w$ for all $g\in U$.
Our aim in this section is to explain how to obtain a formula for~$w$
in terms of $x\in R$, the elements of~$U$, and the values of~$\beta$
(this is Task~3 of the previous section).

Let us start with the case when $U = C_p$ is a cyclic group of order~$p$.
We denote $\sigma$ a generator of~$U$.
If $\beta : U\to R$ is a 1-cocycle of $U$ with values in~$R$,
then by Lemma~3.2~(c)
$$N_U(\beta(\sigma))  = (1 + \sigma + \cdots + \sigma^{p-1})\beta(\sigma) =
0.$$
Now by Lemma~1 of~[3] we have $\beta(\sigma) = (\sigma - 1) w$,
where
$$w = \sum_{k = 1}^{p-1}\, (1 + \sigma + \cdots + \sigma^{k-1})
\bigl(x \sigma^{-k} \beta(\sigma)\bigr) \in R. \eqno (5.1)$$
The right-hand side of~(5.1) is a noncommutative polynomial with
integer coefficients in the variables~$u(x)$ and $u(\beta(\sigma))$, 
where $u\in U$.
By Lemma~3.2~(b) we obtain
for $g = \sigma^i$, where $1\leq i \leq p-1$,
$$\eqalign{
\beta(g) & = (1 + \sigma + \cdots + \sigma^{i-1})\beta(\sigma)\cr
& = (1 + \sigma + \cdots + \sigma^{i-1})(\sigma - 1) w \cr
& = (\sigma^i - 1) w = (g-1)w. \cr
}$$
With Formula~(5.1) we have thus expressed any 1-cocycle as a 1-cobound\-ary
in the case when $U$ is a cyclic group of order~$p$.
Formula~(5.1) is the prototype of formulas we wish to obtain for $w$
in the general case.

To deal with a general finite $p$-group $U$,
we proceed by induction on the order of~$U$.
Fix a normal subgroup $U'$ of~$U$ of index~$p$
and choose $\sigma \in U$ such that its class
$\bar{\sigma}$ in~$U/U'$ generates~$U/U'$.
Following~(1.2), set
$$x' = (1 + \sigma + \cdots + \sigma^{p-1})(x). \eqno (5.2)$$
Then~$N_{U'}(x') = 1$.
We assume that we know how to express any 1-cocycle $\gamma : U'\to R$ 
as the coboundary of an element of~$R$ expressed as a noncommutative 
polynomial with integer coefficients 
in $u'(x')$ and $u'(\gamma(v'))$~($u',v' \in U'$).

In order to pass from $U'$ to~$U$ we make use of a well-known
construction due to Wall~[9].
Let $(C'_*,d')$ be the standard resolution of $\ZZ$ by
free left $\ZZ[U']$-modules.
In particular, $C'_0  = \ZZ[U']$, $C'_1 = \ZZ[U'\times U']$ and the
differential
$d' : C'_1 \to C'_0$ is given for all $g$, $h\in U'$ by
$$d'(g,h) = gh - g.$$
For each $q \geq 0$ we define a chain complex $C_{*,q}$ of free left
$\ZZ[U]$-modules
by setting
$$C_{p,q} = \ZZ[U] \otimes_{\ZZ[U']} C'_p.$$
We define a differential $d_0 : C_{p,q} \to C_{p-1,q}$ by
$d_0  = \id_{\ZZ[U]} \ot d'$.
Observe that
$$C_{0,q} = \ZZ[U] \and C_{1,q} = \ZZ[U\times U']$$
for all $q\geq 0$.
The chain complex $C_{*,q}$ is a free resolution of
$\ZZ[U] \otimes_{\ZZ[U']} \ZZ$,
which can be identified with $\ZZ[U/U']$.
By Lemma~2 and Theorem~1 of~[9] there exist $\ZZ[U]$-linear maps
$$d_k : C_{p,q} \to C_{p+k-1,q-k} \qquad (k\geq 1, q\geq k)$$
such that

(i) when $p=0$, then $d_1 : C_{0,q} = \ZZ[U] \to C_{0,q-1} = \ZZ[U]$ is
given by
$$d_1(\xi) = \xi(1 + \sigma + \cdots + \sigma^{p-1})
\;\;\hbox{if}\; \xi \in C_{0,2i}, \eqno (5.3)$$
$$d_1(\xi) = \xi(\sigma - 1) \;\;\hbox{if}\; \xi \in C_{0,2i-1} \eqno (5.4)$$
(here $i\geq 1$), and

(ii)
$$\sum_{i=0}^k\, d_id_{k-i} = 0. \eqno (5.5)$$

\smallskip\noindent
Define a nonnegatively graded $\ZZ[U]$-module $C_*^{\rm W}$ for all $r\geq
0$ by
$$C_r^{\rm W} = \bigoplus_{p+q = r}\, C_{p,q}.$$
Observe that
$$C_0^{\rm W} = \ZZ[U] \and C_1^{\rm W} = \ZZ[U\times U'] \oplus \ZZ[U].$$
The maps $d^{\rm W} = \sum_{k\geq 0}\, d_k$ define a degree~$-1$
differential on
$C_*^{\rm W}$ and turn it into a resolution of~$\ZZ$ by free left
$\ZZ[U]$-modules.

Let us apply the functor $\Hom_{\ZZ[U]}(-,R)$ to the resolution $(C^{\rm
W}_*,d^{\rm W})$.
Define
$$C^{p,q}  = \Hom_{\ZZ[U]}(C_{p,q},R) =
\Hom_{\ZZ[U]}(\ZZ[U] \otimes_{\ZZ[U']} C'_p,R)
= \Hom_{\ZZ[U']}(C'_p,R)$$
(the last isomorphism follows by adjunction).
In particular,
$$C^{0,q} = R \and C^{1,q} = \Hom(U',R)$$
for all~$q\geq 0$.
The differential $d_0$ turns into a degree~$+1$ differential
$\delta^0 : C^{p,q} \to C^{p+1,q}$.
The maps $d_k$ ($k\geq 1$) turn into maps $\delta^k : C^{p,q} \to
C^{p-k+1,q+k}$,
which by~(5.5) satisfy
$$\sum_{i=0}^k\, \delta^i \delta^{k-i} = 0. \eqno (5.6)$$
Set
$$C^p_{\rm W} = \bigoplus_{i=0}^p\, C^{i,p-i}
\and
\delta_{\rm W} = \sum_{k\geq 0}\, \delta^k .$$
Then $(C^*_{\rm W},\delta_{\rm W})$ is a cochain complex
whose cohomology groups are the groups $H^*(U,R)$.

Any element of $H^1(U,R)$ can be represented by a 1-cocycle
in the cochain complex $(C^*_{\rm W},\delta_{\rm W})$, namely by a couple
$$(\gamma,s) \in C^{1,0} \times C^{0,1} = \Hom(U',R) \times R$$
satisfying
$$\delta^0(\gamma) = 0, \quad \delta^1(\gamma) + \delta^0(s) = 0,
\quad \delta^2(\gamma) + \delta^1(s) = 0. \eqno (5.7)$$
Here $\delta^1(s) = (1 + \sigma + \cdots + \sigma^{p-1}) s$.
A 1-coboundary in the complex $(C^*_{\rm W},\delta_{\rm W})$ is a couple
$(\gamma,s) \in C^{1,0} \times C^{0,1} = \Hom(U',R) \times R$ for which
there exists $w \in C^{0,0} = R$ such that
$$\gamma = \delta^0(w) \and s = \delta^1(w) = (\sigma - 1) w.\eqno (5.8)$$

Let us explain how to find $w\in R$ for a given 1-cocycle~$(\gamma,s)$.
For each $q\geq 0$, $(C^{*,q}, \delta^0)$ is the standard
cochain complex whose cohomology groups are the groups $H^*(U',R)$.
In particular, the kernel of $\delta^0  : R = C^{0,q} \to C^{1,q}$ is~$R^{U'}$.
By the first relation in~(5.7) the element $\gamma \in C^{1,0}$ is a
1-cocycle for the
cochain complex~$(C^{*,0}, \delta^0)$.
By assumption we know how to construct $w_1 \in R$
such that $\gamma = \delta^0(w_1)$ polynomially in terms
of the norm one element~$x'$, the values of~$\gamma$, 
and the elements of~$U'$.

Set $s' = s - \delta^1(w_1) = s - (\sigma - 1) w_1 \in R$.
Then by~(5.6) and by the second relation in~(5.7),
$$\delta^0(s') = \delta^0(s) - \delta^0\delta^1(w_1)
=  \delta^0(s) + \delta^1\delta^0(w_1) = \delta^0(s) + \delta^1(\gamma) = 0.$$
This proves that $s'$ belongs to~$R^{U'}$.

The element $x'' = N_{U'}(x)$ belongs to $R^{U'}$ and we have
$$(1 + \sigma + \cdots + \sigma^{p-1})x'' =
(1 + \sigma + \cdots + \sigma^{p-1}) N_{U'}(x) = N_U(x) = 1.\eqno (5.9)$$
The third relation in~(5.7), together with (5.3) and~(5.6), implies
$$\eqalign{
(1 + \sigma + \cdots + \sigma^{p-1}) s'
& = \delta^1(s')
= \delta^1(s) - \delta^1\delta^1(w_1) \cr
& = \delta^1(s) + \delta^2\delta^0(w_1) \cr
&= \delta^1(s) + \delta^2(\gamma) = 0. \cr
}$$
Since $\sigma^p$ belongs to the subgroup~$U'$, the element $\sigma^p -1$
acts as $0$
on~$R^{U'}$, and $\sigma$ generates a cyclic group of order~$p$
in the automorphism group of the ring~$R^{U'}$.
The element $s'\in R^{U'}$ is of norm zero for this cyclic group.
Using Formula~(5.1),
we obtain an element $w_2 \in R^{U'}$ such that $s' = (\sigma - 1)w_2$,
explicitly in terms of $s'$, of $\sigma$, and of the element $x''$
appearing in~(5.9).

We claim that $w = w_1 + w_2 \in R$ satisfies Equations~(5.8). Indeed,
$$\delta^1(w) = (\sigma - 1) w = (\sigma - 1) w_1 + (\sigma - 1) w_2
= (\sigma - 1) w_1 + s' = s.$$
On the other hand, $\delta^0(w_2) = 0$ since $w_2$ belongs to~$R^{U'}$.
Therefore,
$$\delta^0(w) = \delta^0(w_1) = \gamma.$$
This proves our claim and shows how to construct~$w$
for the cochain complex $(C^*_{\rm W},\delta_{\rm W})$.

In order to express a 1-cocycle $\beta : U\to R$ in the {\it standard}
cochain complex
as a 1-coboundary, we use the comparison lemma between
the resolution $(C_*^{\rm W}, d^{\rm W})$ and the standard resolution
$(C_*,d)$ of~$\ZZ$ by free left $\ZZ[U]$-modules
(see, e.g., Proposition~1.2 in [6, Chapter~V]).

\medskip\goodbreak
\noindent
{\sc Lemma~5.1}.---
{\it There exists a chain map
$$\theta_* : (C^{\rm W}_*,d^{\rm W}) \to (C_*,d)$$
such that
$\theta_0 : C^{\rm W}_0  = \ZZ[U] \to C_0 = \ZZ[U]$
is the identity map and
$$\theta_1 : C^{\rm W}_1  = \ZZ[U\times U']\oplus \ZZ[U] \to C_1 =
\ZZ[U\times U]$$
is the $\ZZ[U]$-linear map whose restriction to the first summand
$\ZZ[U\times U']$
is induced by the natural inclusion of~$U'$ into~$U$,
and the restriction to the second summand $\ZZ[U]$ is defined for all
$g\in U$ by $\theta_1(g) = (g,\sigma) \in C_1$.
}
\medskip

\Pr
The existence of $\theta_*$ follows from the comparison lemma.
Since $C_0^W = C_0 = \ZZ[U]$,
we can take $\theta_0$ to be the identity map.
It now suffices to check that $d\theta_1 = d^{\rm W}$ for the
map $\theta_1$ described in the lemma.
On the summand $\ZZ[U\times U']$ we have
$$d(\theta_1(g,h)) = d(g,h) = gh - g = d_0(g,h) = d^{\rm W}(g,h)$$
for $g\in U$ and $h\in U'$. On the summand $\ZZ[U]$, by~(5.4) we have
$$d(\theta_1(g)) = d(g,\sigma) = g(\sigma - 1)
= g\sigma - g = d_1(g) = d^{\rm W}(g)$$
for $g\in U$.
\hfill\cqfd
\medskip

When we apply the functor $\Hom_{\ZZ[U]}(-,R)$
to $\theta_* : (C^{\rm W}_*,d^{\rm W}) \to (C_*,d)$, we obtain
a cochain map
$$\theta^* : C^*(U,R)  =  \Hom_{\ZZ[U]}(C_*,R)\to
C^*_{\rm W} = \Hom_{\ZZ[U]}(C^{\rm W}_*,R)$$
inducing an isomorphism in cohomology.
The cochain complex $(C^*(U,R), \delta)$ is the standard cochain complex
computing~$H^*(U,R)$. Now, let $\beta: U\to R$ be a standard 1-cocycle.
This is an element of~$C^1(U,R)$ such that $\delta(\beta) = 0$.
Consider its image $\theta^1(\beta)\in C^1_{\rm W}$. It is a 1-cocycle
in~$(C^*_{\rm W},\delta_{\rm W})$.
By our investigation above we know how to construct $w\in R$
such that $\theta^1(\beta) = \delta_{\rm W}(w)$. We claim the following.

\medskip\goodbreak
\noindent
{\sc Lemma~5.2}.---
{\it We have $\beta = \delta(w)$.
}
\medskip

\Pr
By construction of~$\theta_1$ we have $d\theta_1 = d^{\rm W}$.
Dualizing, we obtain $\theta^1\delta = \delta_{\rm W}$. Therefore,
$$\theta^1(\delta(w)) = \delta_{\rm W}(w) = \theta^1(\beta).$$
To conclude it suffices to check that $\theta^1$ is injective.
Using the string of natural isomorphisms
$$\eqalign{
C^1_{\rm W} & = \Hom_{\ZZ[U]}(C^{\rm W}_1,R) \cr
& = \Hom_{\ZZ[U]}(\ZZ[U\times U'],R) \oplus \Hom_{\ZZ[U]}(\ZZ[U],R)\cr
& = \Hom(U',R) \oplus R\cr
}$$
and Lemma~5.1, we easily see that the image $\theta^1(\beta)$
of any standard 1-cocycle $\beta \in \Hom(U,R)$ is given by
$$\theta^1(\beta) = (\beta',\beta(\sigma)) \in \Hom(U',R) \oplus R
= C^1_{\rm W},$$
where $\beta'$ is the restriction of~$\beta$ to~$U'$
and $\beta(\sigma)$ is its value on~$\sigma$.
If $\theta^1(\beta) = 0$, then
the restriction of~$\beta$ to~$U'$ is zero and $\beta(\sigma) = 0$.
It follows from Lemma~3.2~(b)
that $\beta$ vanishes on all powers of~$\sigma$. The cocycle condition~(3.1)
then implies that $\beta$ vanishes on all elements of~$U$.
This proves the injectivity of~$\theta_1$.
\hfill\cqfd
\medskip

Summing up, we thus have obtained an inductive way (starting from cyclic
groups) to
express any 1-cocycle of a finite $p$-group (with values in a ring~$R$)
as the coboundary of an element $w\in R$, 
polynomially in terms of~$x$,
the values of the $1$-cocycle, and the elements of the group.
This is a vast generalization of~[3, Lemma~1].

\medskip\goodbreak
\noindent
{\sc Example~5.3}. Let $U$ be an elementary abelian group generated by two
generators
$\sigma_1$ and $\sigma_2$ of order two and acting on a ring~$R$.
Let $U'$ be the subgroup generated by $\sigma_1$.
We assume the existence of an element $x\in R$ such that $N_U(x) = 1$.
The elements
$$x' = (1 + \sigma_2)(x) \and x'' = (1 + \sigma_1)(x)$$
are of norm one for $U'$ and $U/U'$, respectively.
Observe that $\sigma_2(x') = x'$.

A 1-cocycle in the complex $(C^*_{\rm W}, \delta_{\rm W})$
corresponding to this situation
is a couple $(\gamma,s) \in  \Hom(U',R) \times R$ satisfying
Equations~(5.7). In particular, $\gamma : U'\to R$
is a 1-cocycle for the subgroup~$U'$.
Set~$r = \gamma(\sigma_1)\in R$. Then Equations~(5.7) are equivalent to
the following three equations:
$$(1 + \sigma_1)(r) = 0, \quad (\sigma_2 - 1)(r) + (\sigma_1 - 1)(s) = 0,
\quad (1 + \sigma_2)(s) = 0.$$
(In this example as in any case when $U$ is a semidirect product of $U'$
and $U/U'$,
the map $\delta^2$ in~(5.7) vanishes.)
By (5.1) we have
$r = (\sigma_1 - 1) w_1$, where
$$w_1 = x' \sigma_1(r) = (1 + \sigma_2)(x)\, \sigma_1(r).$$
Consequently,
$$s' = s - (\sigma_2 - 1) w_1 = s - x' (\sigma_2\sigma_1)(r) + x' \sigma_1(r)
 = s - x'(\sigma_1(\sigma_2 -1)(r)).$$
Following the procedure above,
we have $s' = (\sigma_2 - 1) w_2$,
where by (5.1)
$$\eqalign{
w_2
& = x''\sigma_2(s') \cr
& = x''\sigma_2 \bigl( s - x'(\sigma_1(\sigma_2 -1)(r)) \bigr) \cr
& = x''\sigma_2(s) - x'' x'(\sigma_2\sigma_1(\sigma_2 -1)(r)) \cr
& = x''\sigma_2(s) + x'' x'(\sigma_1(\sigma_2 -1)(r)) \cr
& = (1 + \sigma_1)(x)\cdot \sigma_2(s)
+ (1 + \sigma_1)(x)\cdot (1 + \sigma_2)(x)\cdot (\sigma_1(\sigma_2 -1)(r)).\cr
}$$
Therefore, if we set
$$\eqalign{
w = & \ w_1 + w_2 \cr
= & \; (1 + \sigma_2)(x)\cdot \sigma_1(r)
+ (1 + \sigma_1)(x)\cdot \sigma_2(s) \cr
& + (1 + \sigma_1)(x)\cdot (1 + \sigma_2)(x)\cdot (\sigma_1(\sigma_2
-1)(r)), \cr
} \eqno (5.10)$$
we obtain
$\gamma(\sigma) = r = (\sigma_1 - 1) w$ and $s = (\sigma_2 - 1) w$.

\bigskip\goodbreak
\noindent
{\bf 6.~The quaternion $2$-groups}.
\medskip

The smallest nonabelian $p$-groups
are the two nonabelian groups of order~$8$,
namely the quaternion group~$Q_8$, which has a unique
elementary abelian subgroup of order~2, and
the dihedral group~$D_8$, which has two nonconjugate
maximal elementary abelian subgroups of order~4.
Both $Q_8$ and $D_8$ are extraspecial groups.

In this section we apply the method of Section~4 in order to
solve the problem for $Q_8$ and more generally
for the generalized quaternion groups $Q_{2^{n+2}}$ ($n\geq 1$).

The group $G = Q_{2^{n+2}}$ of order $2^{n+2}$ (with $n\geq 1$) has a
presentation
with two generators $\sigma$, $\tau$ and the relations
$$\sigma^{2^{n+1}} = 1, \quad \tau\sigma = \sigma^{-1}\tau,
\quad \tau^2 = \sigma^{2^n} . \eqno (6.1)$$
Any element of the group can be written as $\sigma^i \tau^j$,
where $i = 0, 1, \ldots, 2^{n+1} -1$ and~$j=0,1$.
We take $U$ to be the cyclic group generated by~$\sigma$.
The quotient group $G/U$ is cyclic of order~$2$ and generated by the class
of~$\tau$.

The group $Q_{2^{n+2}}$ has a unique elementary abelian subgroup, which is the
group of order~2 generated by the central element~$\tau^2 = \sigma^{2^n}$.

We follow the method presented in Section~4.
Let us first perform Task~1.

\medskip\goodbreak
\noindent
{\sc Lemma~6.1}.---
{\it The values $b(\sigma)$ and $b(\tau) \in R$ satisfy the system of three
equations
$$
\left\{\matrix{
N_U(b(\sigma)) & = 0, \cr
\noalign{\smallskip}
(\sigma - 1) b(\tau) + (1 + \sigma\tau) b(\sigma) & = 0, \cr
\noalign{\smallskip}
(1 + \tau) b(\tau) - (1 + \sigma + \cdots + \sigma^{2^n -1}) b(\sigma) & =
1. \cr
}\right.
$$
}

\Pr
Since the restriction of~$b$ to $U$ is a 1-cocycle, the first equation follows
from Lemma~3.2~(c).
Applying Lemma~3.2~(e) to $\beta = b - \alpha$, we obtain
$$ (\sigma - 1)(b(\tau) - \varphi) + (1 + \sigma\tau) b(\sigma) = 0.$$
We derive the second equation of the lemma by recalling that $\varphi$
is $U$-invariant.
In order to prove the third equation, we use the third relation in~(6.1).
Since $b$ is a 1-cocycle when restricted to~$U$, we have
$$b(\sigma^{2^n}) = (1+ \sigma + \cdots + \sigma^{2^n -1}) b(\sigma) \eqno
(6.2)$$
by Lemma~3.2~(b).
On the other hand, we have
$$b(\tau^2) = (1 + \tau) (b(\tau) - \varphi)
= (1 + \tau) b(\tau) - N_{G/U}( \varphi )
= (1 + \tau) b(\tau) - 1.
\eqno (6.3)$$
The third equation then follows from $\tau^2 = \sigma^{2^n}$ and~(6.2--6.3).
\hfill\cqfd
\goodbreak\medskip

To solve Task~2 , we need an element $x$ of~$R$ such that
$$N_U(x)  = (1 + \sigma + \cdots + \sigma^{2^{n+1} -1})(x) = 1.$$

\medskip\goodbreak
\noindent
{\sc Lemma~6.2}.---
{\it The elements
$$b(\sigma) = (1 - \sigma\tau)(x) \and
b(\tau) = (1+\sigma + \cdots + \sigma^{2^n -1})(x)$$
of~$R$ are solutions of the system of equations of~Lemma~6.1.
}
\medskip

\Pr
For the first equation we have
$$\eqalign{
N_U(b(\sigma))
& = N_U(1 - \sigma\tau)(x) = N_U(x) - N_U\tau(x) \cr
& = N_U(x) - \tau N_U(x) = (1 - \tau)(1) = 0.\cr
}$$
We check the second equation:
$$\eqalign{
(\sigma - 1) b(\tau) + (1 + \sigma\tau) b(\sigma)
& = (\sigma - 1) (1+\sigma + \cdots + \sigma^{2^n -1})(x) \cr
& \qquad {} + (1 + \sigma\tau) (1 - \sigma\tau)(x) \cr
& = (\sigma^{2^n} - 1 + 1 - (\sigma\tau)^2)(x) = 0. \cr
}$$
For the third equation we have
$$\eqalign{
(1 + \tau) b(\tau) & - (1 + \sigma + \cdots + \sigma^{2^n -1}) b(\sigma) \cr
& = \bigl( (1 + \tau) (1+\sigma + \cdots + \sigma^{2^n -1})
- (1 + \sigma + \cdots + \sigma^{2^n -1}) (1 - \sigma\tau) \bigr)(x) \cr
& = \tau (1+ \sigma + \cdots + \sigma^{2^{n+1} -1})(x)  \cr
& = \tau(N_U(x)) = \tau(1) = 1.\cr
}$$
\line{\hfill\cqfd}
\goodbreak\medskip

We now complete Task~3, which is to find an explicit $w\in R$ such that
$b(g) = (g-1)w$ for $g\in U$.
By Lemmas~6.1--6.2 we have $N_U(b(\sigma)) = 0$ for
$b(\sigma) = (1 - \sigma\tau)(x)$.
Since $U$ is cyclic, we may apply~[3, Lemma~1].
We then obtain $b(\sigma) = (\sigma - 1)w$, where
$$\eqalign{
w
& = \sum_{k=1}^{2^{n+1}-1}\,
(1 + \sigma + \cdots + \sigma^{k-1})\bigl( x\, \sigma^{-k}b(\sigma) \bigr) \cr
& = \sum_{k=1}^{2^{n+1}-1}\,
(1 + \sigma + \cdots + \sigma^{k-1})\bigl( x\, \sigma^{-k}(1 -
\sigma\tau)(x) \bigr). \cr
} \eqno (6.5)$$
Observe that $w$ is a noncommutative polynomial with
$2^{n+1}(2^{n+1}-1)$ monomials of degree $\leq 2$ in terms of~$x$.

As a consequence of Proposition~4.4,
the element $a = b(\tau) + (1 - \tau) w \in R$
is $U$-invariant and we have $N_{G/U}(a) = 1$.
Therefore, $N_G(ax) = 1$ for~$G = Q_{2^{n+2}}$.
It can be checked that $y = ax$ is a polynomial in the variables
$g(x)$ ($g\in G$) with
$2^n(1 + 4 (2^{n+1}-1))$ monomials of degree~$\leq 3$.

For the special case when $n=1$ and $G = Q_8$ is the quaternion
group of order~8, we obtain the following element $y \in R$
satisfying $N_{Q_8}(y) = 1$:
$$\eqalign{
y & =  x^2 + \sigma(x)x \cr
&\quad + x\sigma(x)x + x\sigma^2(x)x + x\sigma^3(x)x \cr
&\quad - x\tau(x) x - x(\sigma^2\tau)(x) x - x (\sigma^3\tau)(x) x \cr
&\quad + \sigma(x)\sigma^2(x)x + \sigma(x)\sigma^3(x)x \cr
&\quad - \sigma(x)\tau(x) x - \sigma(x)(\sigma^3\tau)(x)x \cr
&\quad + \sigma^2(x)\sigma^3(x)x - \sigma^2(x)\tau(x)x \cr
&\quad + \tau(x)x^2 + \tau(x)\sigma^2(x)x + \tau(x)\sigma^3(x)x \cr
&\quad  -  \tau(x)(\sigma\tau)(x)x - \tau(x)(\sigma^2\tau)(x)x
	 - \tau(x)(\sigma^3\tau)(x) x \cr
&\quad + (\sigma^2\tau)(x)\sigma^2(x) x - (\sigma^2\tau)(x)(\sigma\tau)(x)
x  \cr
&\quad  + (\sigma^3\tau)(x) \sigma^2(x) x + (\sigma^3\tau)(x) \sigma^3(x) x \cr
&\quad - (\sigma^3\tau)(x)(\sigma\tau)(x) x
- (\sigma^3\tau)(x)(\sigma^2\tau)(x) x .\cr
} \eqno (6.6)$$
The right-hand side of (6.6) contains~26 monomials of degree~$\leq 3$
in terms of~$x$. 
If we wish to express $y$ in terms of a norm one element~$x_E$
for the elementary abelian subgroup~$E$ of~$Q_8$ generated by~$\sigma^2$,
it suffices by~(0.1) to replace $x$ in~(6.6) by the polynomial
$x_E\sigma(x_E)x_E + x_E\sigma(x_E) - x_E^2\sigma(x_E)$.
We thus obtain a formula for~$Q_8$ with 
$666$ $( = 2 \cdot 3^2 + 24 \cdot 3^3)$ monomials of degree~$\leq 9$.

\medskip\goodbreak
\noindent
{\sc Remark~6.3}. Observe that the group $Q_8$ is extraspecial, and
if $G = Q_{2^{n+2}}$ ($n\geq 2$), then $\cF_G = \{C_4, Q_8, D_8 \}$.

\goodbreak
\bigskip\goodbreak
\noindent
{\bf 7.~The dihedral $2$-groups}
\medskip

The dihedral group $G = D_{2^{n+1}}$ of order $2^{n+1}$ (where $n\geq 2$)
has a presentation with two generators $\sigma$, $\tau$
and the relations
$$\tau\sigma = \sigma^{-1}\tau \and \tau^2 = \sigma^{2^n} = 1 . \eqno (7.1)$$
Any element of the group
can be written uniquely as $\sigma^i \tau^j$, where $i = 0, 1,\ldots, 2^n
-1$ and~$j=0,1$.
Let $U$ be the normal subgroup generated by $\sigma^2$ and~$\tau$.
The quotient group $G/U$ is the cyclic group of order~2
generated by the class of~$\sigma$.
Note that $U$ is a dihedral group of order~$2^n$ if $n\geq 3$ and
an elementary abelian group if~$n=2$. It contains the elementary
abelian subgroup $U_1$ generated by $u$ and~$\tau$, where $u =
\sigma^{2^{n-1}}$
is the unique non-trivial central element of $D_{2^{n+1}}$.

Let $x$ be an element of~$R$ such that $N_{U}(x)  = 1$.
We denote $H$ the cyclic group generated by~$\sigma$ (of order~$2^n$).
We follow the method presented in Section~4.
Let us first perform Task~1.

\medskip\goodbreak
\noindent
{\sc Lemma~7.1}.---
{\it The values $b(\sigma)$ and $b(\tau) \in R$ satisfy the system of three
equations
$$
\left\{\matrix{
(1+\tau) b(\tau) & = & 0, \cr
\noalign{\smallskip}
N_H(b(\sigma))  & = & 2^{n-1}, \cr
\noalign{\smallskip}
(\sigma - 1) b(\tau) + (1 + \sigma\tau) b(\sigma) & = &1. \cr
}\right.
$$
}
\medskip

\Pr
Since $b: U\to R$ is a 1-cocycle, and $\tau$ and $\sigma^2$ belong to~$U$,
we have
$$(1 + \tau)b(\tau) = 0 \and
(1 + \sigma^2 + \cdots + \sigma^{2^n-2})b(\sigma^2) = b(\sigma^{2^n}) =
b(1) = 0.
 \eqno (7.2)$$
This proves the first equation.
We have
$$(\alpha - b)(\sigma^2) = (1 + \sigma)(\alpha - b)(\sigma)$$
by Lemma~3.2~(b); hence
$$b(\sigma^2) = (1+ \sigma)( b(\sigma) - \varphi)
= (1+ \sigma)b(\sigma) - N_{G/U}(\varphi)
= (1+ \sigma)b(\sigma) - 1. \eqno (7.3)$$
The second relation in~(7.2) and Relation (7.3) imply
$$\eqalign{
N_H(b(\sigma))
& = (1 + \sigma^2 + \cdots + \sigma^{2^n-2}) (1+ \sigma)b(\sigma) \cr
& = (1 + \sigma^2 + \cdots + \sigma^{2^n-2}) b(\sigma^2)
+ (1 + \sigma^2 + \cdots + \sigma^{2^n-2})(1) \cr
& = 2^{n-1}.\cr
}$$
The second equation of the lemma is thus proved.
Applying Lemma~3.2~(e) to $\beta = b - \alpha$, we obtain
$$ (\sigma - 1)b(\tau)  + (1 + \sigma\tau) (b(\sigma) - \varphi) = 0,$$
which implies
$$(\sigma - 1) b(\tau) + (1 + \sigma\tau) b(\sigma)  = (1 + \sigma\tau)\varphi
 = N_{G/U}(\varphi) = 1.$$
This proves the last equation.
\hfill\cqfd
\medskip

Let $U_2$ be the elementary abelian subgroup of~$G = D_{2^{n+1}}$
generated by $u$ and~$\sigma\tau$.
(The subgroups $U_1$ and $U_2$ are not conjugate in~$G$.)
Let $x_2$ be an element of~$R$ satisfying
$$N_{U_2}(x_2)  = (1 + u + \sigma\tau + u\sigma\tau)(x_2)
= (1 + \sigma\tau) (1+ u)(x_2) = 1. \eqno (7.4)$$

\medskip\goodbreak
\noindent
{\sc Lemma~7.2}.---
{\it The elements
$$b(\sigma) = (\sigma + u\tau)(x_2) \and
b(\tau) = (\tau-1)(1+u)(x_2)$$
of~$R$ are solutions of the system of equations of~Lemma~7.1.
}
\medskip

\Pr
The first equation is clearly satisfied. For the second one, we have
$$\eqalign{
N_H(b(\sigma))
& = (N_H\sigma)(x_2) + (N_Hu\tau)(x_2) = N_H(1 + \tau)(x_2) \cr
& = N_G(x_2) = [G:U_2] \, N_{U_2}(x_2) = 2^{n-1}.\cr
}$$
We now check the third equation. Using (7.4) and the identities
$$(\sigma-1)(\tau-1) = (1 + \sigma\tau)(1 - \sigma)
\and
N_{U_2} = (1 + \sigma\tau)(1+u)$$
in $\ZZ[G]$, we obtain
$$\eqalign{
(\sigma - 1) b(\tau) + (1 + \sigma\tau) b(\sigma) -1
& = (1 + \sigma\tau)
\bigl( (1 - \sigma)(1+u) + (\sigma + u\tau) - (1+u)\bigr)(x_2) \cr
& = (1 + \sigma\tau)(\tau -\sigma)u(x_2) \cr
& = (1 + \sigma\tau)(\sigma\tau -1)\sigma u(x_2) \cr
& = ((\sigma\tau)^2 -1)\sigma u(x_2) = 0. \cr
}$$
\line{\hfill\cqfd}
\goodbreak\medskip

Proceeding as in Section~5, we can find $w\in R$ such that
$$b(\sigma^2) = (\sigma^2 -1) w \and b(\tau) = (\tau -1)w.$$
The element $w$ can be expressed (as a noncommutative polynomial
with integer coefficients) in terms of the norm one
element~$x$, the elements of~$U$, and the values $b(\sigma)$, $b(\tau)$
given in Lemma~7.2.
Observe that here we need both $x$ and $x_2\in R$, which is not surprising
since $G$ has two nonconjugate maximal elementary abelian subgroups.
As a consequence of Proposition~4.4, the element
$$a  = b(\sigma) + (1 - \sigma)(w)$$
is $U$-invariant and $N_{G/U}(a) = 1$.
Hence, $N_G(y) = 1$ for $y = ax$ by Proposition~4.1.

If $G = D_8$, then $U = U_1$ is elementary abelian of order~$4$,
and we can use Example~5.3.
Setting $\sigma_1 = \tau$, $\sigma_2 = \sigma^2$, $r = b(\tau)$,
and $s = b(\sigma^2)$ in Formula~(5.10),
we obtain an explicit $w$ with 48 monomials of degree at most~3.
Hence for $G = D_8$ we have a norm one element~$y$
with 98 monomials of degree~$\leq 4$.

\medskip\goodbreak
\noindent
{\sc Remarks~7.3}. (a) The group $D_8$ is extraspecial and,
if $G = D_{2^{n+1}}$ ($n\geq 2$), then $\cF_G = \{C_4, D_8\}$.

(b) By Sections 2, 6, 7, we have solved the problem for all $2$-groups $G$
such that $\cF_G = \{C_4, Q_8, D_8\}$, in particular for 
all metacyclic $2$-groups. Note that by [4, Theorem~5.1]
any $2$-group every subgroup of which is generated by two elements is metacyclic.

\bigskip\goodbreak
\noindent
{\bf 8.~A nonabelian group of order~$27$}
\medskip

Let $p$ be an odd prime number and
$G_{p^3}$ be the group generated by $\sigma$, $\tau$ and the relations
$$\sigma^{p^2} = \tau^p = 1\and  \tau\sigma = \sigma^{p+1}\tau . \eqno (8.1)$$
This is the only nonabelian group of order $p^3$ containing a
cyclic subgroup of index~$p$. The center $Z$ of~$G_{p^3}$ is the cyclic group
generated by $\sigma^p$, and $G_{p^3}/Z$ is elementary abelian
of order~$p^2$. Therefore, $G_{p^3}$ is extraspecial.

Let $U$ be the elementary abelian subgroup of~$G_{p^3}$ generated by
$\sigma^p$ and~$\tau$; it is the unique maximal elementary abelian subgroup
of~$G_{p^3}$.
The quotient group $G_{p^3}/U$ is the cyclic group of order~$p$
generated by the class of~$\sigma$.

Let $x$ be an element of~$R$ such that $N_{U}(x)  = 1$.
We denote $H$ the cyclic group of order~$p^2$ generated by~$\sigma$.
Following the method of Section~4, we undertake Task~1.

\medskip\goodbreak
\noindent
{\sc Lemma~8.1}.---
{\it The values $b(\sigma)$ and $b(\tau) \in R$ satisfy the system of three
equations
$$
\left\{\matrix{
(1+\tau + \cdots + \tau^{p-1}) b(\tau) & = & 0, \cr
\noalign{\smallskip}
N_H(b(\sigma))  & = & p, \cr
\noalign{\smallskip}
(\sigma^{p+1} - 1) b(\tau)
+ (1 + \sigma + \ldots + \sigma^{p-1} + \sigma^p - \tau) b(\sigma) & = &1. \cr
}\right.
$$
}
\medskip

\Pr
Since $b: U\to R$ is a 1-cocycle, and $\tau$ and $\sigma^p$ belong to~$U$,
we have
$$(1 + \tau + \cdots + \tau^{p-1})b(\tau) = 0 \and
(1 + \sigma^p + \cdots + \sigma^{(p-1)p})b(\sigma^p) = 0.
 \eqno (8.2)$$
This proves the first equation.
By Lemma~3.2~(b) we have
$$(\alpha - b)(\sigma^p) = (1 + \sigma + \cdots + \sigma^{p-1})(\alpha -
b)(\sigma),$$
which implies
$$\eqalign{
b(\sigma^p) & = (1+ \sigma + \cdots + \sigma^{p-1})( b(\sigma) - \varphi) \cr
& = (1+ \sigma + \cdots + \sigma^{p-1})b(\sigma) - N_{G_{p^3}/U}(\varphi) \cr
& = (1+ \sigma + \cdots + \sigma^{p-1})b(\sigma) - 1. \cr
}\eqno (8.3) $$
The second relation in~(8.2), together with Relation (8.3), implies
$$\eqalign{
N_H(b(\sigma))
& = (1 + \sigma^p + \cdots + \sigma^{(p-1)p})
(1+ \sigma + \cdots + \sigma^{p-1})b(\sigma) \cr
& = (1 + \sigma^p + \cdots + \sigma^{(p-1)p})  b(\sigma^p)
+ (1 + \sigma^p + \cdots + \sigma^{(p-1)p})(1) \cr
& = p.\cr
}$$
This proves the second equation.
To prove the last one, we first compute~$b(\sigma^{p+1})$. We have
$$b(\sigma^{p+1}) - \varphi = b(\sigma^p) + \sigma^p (b(\sigma) - \varphi),$$
hence
$$\eqalign{
b(\sigma^{p+1}) & = b(\sigma^p) + \sigma^p b(\sigma) - (\sigma^p -1)\varphi \cr
& = (1+ \sigma + \cdots + \sigma^{p-1} + \sigma^p)b(\sigma) - 1\cr
}$$
in view of~(8.3) and the $\sigma^p$-invariance of~$\varphi$.
Applying the cocycle condition to the third relation in~(8.1),
we obtain
$$\eqalign{
b(\tau) + \tau (b(\sigma) - \varphi) & = b(\tau\sigma) - \varphi
= b(\sigma^{p+1}\tau) - \varphi \cr
& = b(\sigma^{p+1}) - \varphi + \sigma^{p+1} b(\tau)\cr
& = (1+ \sigma + \cdots + \sigma^{p-1} + \sigma^p)b(\sigma) - 1 - \varphi +
\sigma^{p+1}
b(\tau).\cr }$$
This, together with the $\tau$-invariance of~$\varphi$, proves the third
equation
of the lemma.
\line{\hfill\cqfd}
\goodbreak\medskip

We will solve the system of equations of Lemma~8.1 when $p = 3$, i.e.,
for the group $G_{27}$ of order~27,
generated by $\sigma$, $\tau$ and the relations
$$\sigma^{9} = \tau^3 = 1\and  \tau\sigma = \sigma^{4}\tau . \eqno (8.4)$$

The elementary abelian subgroup~$U$ considered above is generated
by $\sigma^3$ and~$\tau$.
We assume the existence of $x\in R$ such that $N_U(x) = 1$.
The center $Z$ of~$G_{27}$ is the cyclic group generated by~$\sigma^3$;
it is contained in~$U$. Therefore, if we set
$x_0 = (1+ \tau + \tau^2)(x) \in R$, we have $N_Z(x_0) = 1$.

Consider the cyclic group $H'$ of order $9$ generated by $\sigma\tau$.
We have $(\sigma\tau)^3 = \sigma^3$. Hence $H'$ contains $Z$ as a
subgroup of index~$3$.
By [3, Corollary~1] we obtain an element $x'\in R$ such that $N_{H'}(x')
= 1$.
To have an explicit formula for~$x'$,
replace $\sigma$ by $\sigma\tau$,
$x_E$ by $x_0 = (1+ \tau + \tau^2)(x)$, and $x_G$ by~$x'$ in Formula~(0.2)
of the introduction.
The element $x' \in R$ is used in the next result.

\medskip\goodbreak
\noindent
{\sc Lemma~8.2}.---
{\it The elements $b(\sigma) = (1 + \sigma\tau + (\sigma\tau)^2 )(x')$ and
$$b(\tau) = (\tau - 1)[\sigma^6 - \sigma(1 + \sigma + \sigma^2 +
\sigma^4)\tau](x')$$
of~$R$ are solutions of the system of equations of~Lemma~8.1.
}
\medskip

\Pr
Set $A = \sigma^6 - \sigma(1 + \sigma + \sigma^2 + \sigma^4)\tau \in
\ZZ[G_{27}]$.
Then $b(\tau) = (\tau - 1)A(x')$.
The first equation in Lemma~8.1 is satisfied because
$$(1 + \tau + \tau^2)b(\tau) = (1 + \tau + \tau^2)(\tau - 1)A(x') = 0.$$
For the second equation we have
$$\eqalign{
N_H(b(\sigma))
& = (1 + \sigma + \sigma^2)(1 + \sigma^3 + \sigma^6)
(1 + \sigma\tau + (\sigma\tau)^2 )(x')\cr
& = (1 + \sigma + \sigma^2)(1 + (\sigma\tau)^3 + (\sigma\tau)^6)
(1 + \sigma\tau + (\sigma\tau)^2 )(x')\cr
& = (1 + \sigma + \sigma^2)N_{H'}(x') = (1 + \sigma + \sigma^2)(1) = 3.\cr
}$$
The following identity in the group ring~$\ZZ[G_{27}]$
can be checked directly:
$$(\sigma^4 - 1) (\tau - 1) A
+ (1 + \sigma + \sigma^2 + \sigma^3 - \tau)
(1 + \sigma\tau + (\sigma\tau)^2 )
= N_{H'}. \eqno (8.5)$$
(This identity was found using a computer.)
Applying both sides of~(8.5) to~$x'$, we obtain the third equation in
Lemma~8.1.
\hfill\cqfd
\medskip

Proceeding as in Example~5.3, we can find $w\in R$ such that
$$b(\sigma^3) = (\sigma^3 -1) w \and b(\tau) = (\tau -1)w.$$
The element $w$ can be expressed (as a noncommutative polynomial
with integer coefficients) in terms of the norm one element~$x$,
the elements of~$U$, and the values $b(\sigma)$, $b(\tau)$
given in Lemma~8.2.
As a consequence of Proposition~4.4, the element
$a  = b(\sigma) + (1 - \sigma)(w)$ is $U$-invariant and $N_{G_{27}/U}(a) = 1$.
Therefore, $N_{G_{27}}(y) = 1$ for $y = ax$ or $y = xa$ by Proposition~4.1.

As a consequence of Section~2, we have solved the problem
for all $3$-groups $G$ such that $\cF_G = \{ C_9, G_{27} \}$.

\bigskip\goodbreak
\noindent
{\sc Acknowledgements}.
We thank Patrick Dehornoy and Ron Holzman for useful discussions on the
best way
to write up the results of Section~2.

\vfill\eject
\bigskip\bigskip\goodbreak
\centerline{\bf References}
\bigskip\bigskip

\noindent
[1] {\sc E.~Aljadeff},
{\it On the surjectivity of some trace maps},
Israel J.~Math.\ 86 (1994), 221--232.
\smallskip

\noindent
[2] {\sc E.~Aljadeff, Y.~Ginosar},
{\it Induction from elementary abelian subgroups},
J.~of Algebra 179 (1996), 599--606.
\smallskip

\noindent
[3] {\sc E.~Aljadeff, C.~Kassel},
{\it Explicit norm one elements for ring actions of finite abelian groups},
Israel J.~Math.\ 129 (2002), 99--108.
\smallskip

\noindent
[4] {\sc N.~Blackburn},
{\it Generalizations of certain elementary theorems on $p$-groups},
Proc.\ London Math.\ Soc.\ 11 (1961), 1--22.
\smallskip

\noindent
[5] {\sc J.\ F.~Carlson, J.~Th\'evenaz},
{\it Torsion endo-trivial modules},
Algebr.\ Represent.\ Theory 3 (2000), 303--335.
\smallskip

\noindent
[6] {\sc H.~Cartan, S.~Eilenberg},
{\it Homological algebra},
Princeton University Press, Princeton, 1956.
\smallskip

\noindent
[7] {\sc L. G. Chouinard},
{\it Projectivity and relative projectivity over group rings},
J.~Pure Appl.\ Algebra 7 (1976), 287--302.
\smallskip

\noindent
[8] {\sc D.~Gorenstein},
{\it Finite groups}, 2nd edition,
Chelsea Publ.\ Co., New York, 1980.
\smallskip

\noindent
[9] {\sc C. T. C. Wall},
{\it Resolutions for extensions of groups},
Proc.\ Cambridge Phil.\ Soc.\ 57 (1961), 251--255.
\smallskip

\vskip 30pt

\line{Eli Aljadeff\hfill}
\line{Department of Mathematics \hfill}
\line{Technion - Israel Institute of Technology \hfill}
\line{32000 Haifa, Israel \hfill}
\line{E-mail:  {\refttfont aljadeff@techunix.technion.ac.il}\hfill}
\line{Fax: +972-4-832-4654 \hfill}
\medskip\medskip

\line{Christian Kassel \hfill}
\line{Institut de Recherche Math\'ematique Avanc\'ee \hfill}
\line{CNRS - Universit\'e  Louis Pasteur\hfill}
\line{7 rue Ren\'e Descartes \hfill}
\line{67084  Strasbourg Cedex, France \hfill}
\line{E-mail:  {\refttfont kassel@math.u-strasbg.fr}\hfill}
\line{Fax: +33 (0)3 90 24 03 28 \hfill}
\line{http:/\hskip
-2.5pt/www-irma.u-strasbg.fr/\raise-4pt\hbox{\~{}}kassel/\hfill}

\bye